\documentclass[11pt]{amsart}
\usepackage[latin9]{inputenc}
\usepackage{color}
\usepackage{cancel}
\usepackage{amsbsy}
\usepackage{amstext}
\usepackage{amsthm}
\usepackage{amssymb}
\PassOptionsToPackage{normalem}{ulem}
\usepackage{ulem}

\makeatletter
\numberwithin{equation}{section}
\numberwithin{figure}{section}

\usepackage{tikz-cd}
\usepackage{fancyhdr}
\usepackage{latexsym}
\usepackage{amscd}\usepackage{amsfonts}\usepackage{pstricks}
\usepackage{young}
\usepackage{enumitem}
\usepackage{amsthm}
\usepackage{mathrsfs}

\usepackage{ifthen}
\usepackage{longtable}
\usepackage{todonotes}
\usepackage{marginnote}

\renewcommand{\subsection}[1]{\vspace{3mm}\refstepcounter{subsection}\noindent{\bf \thesubsection. #1.} }

\renewcommand{\subsubsection}[1]{\vspace{3mm}\refstepcounter{subsubsection}\noindent{\bf \thesubsubsection. #1.} }

\numberwithin{equation}{section}

\newtheorem{theorem}{Theorem}
\newtheorem{lemma}[theorem]{Lemma}
\newtheorem{corollary}[theorem]{Corollary}
\newtheorem{proposition}[theorem]{Proposition}
\newtheorem{conjecture}{Conjecture}

\theoremstyle{definition}

\newtheorem{definition}[theorem]{Definition}

\theoremstyle{remark*}
\newtheorem*{remark*}{Remark}
\newtheorem{example}{Example}
\newtheorem*{example*}{Example}
\newtheorem{notation}{Notation}

\def\CC{\mathbb C}
\def\PP{\mathbb P}

\def\ord{\operatorname{ord}}

\DeclareMathOperator{\Res}{Res}
\DeclareMathOperator{\GL}{GL}

\def\min{\mathop{\mathrm{min}}}
\def\CC{\mathbb C}

\def\ZZ{\mathbb Z}
\def\PP{\mathbb P}
\def\K{K}

\def\cal{\mathcal }
 \def\ord{\text{ord}}
 
\def\p{\mathbf p}
\def\q{\mathbf q}

\def\gen{\mathfrak g}

\let\a\alpha

\newcommand{\bfi}{\mathbf{i}}
\newcommand{\bfj}{\mathbf{j}}
\newcommand{\bfx}{\mathbf{x}}
\newcommand{\bfu}{\mathbf{u}}
\newcommand{\bfv}{\mathbf{v}}

\newcommand{\bft}{\mathbf{t}}

\makeatother

\begin{document}
\title[ Vojta's abc Conjecture  for algebraic tori and applications over Function Fields ]{Vojta's abc Conjecture  for algebraic tori \\ and applications over Function Fields }
\author{Ji Guo}
\address{School of Mathematics  and Statistics\\HNP-LAMA \\ Central South University \\ Changsha  410075 \\ China} 
\email{221250@csu.edu.cn}

\author{Khoa D. Nguyen}
\address{Department of Mathematics and Statistics, University of Calgary, AB T2N 1N4, Canada}
\email{dangkhoa.nguyen@ucalgary.ca}

\author{Chia-Liang Sun}
\address{Institute of Information Science \\Academia Sinica\\128 Academia Road, Section 2\\ Nankang, Taipei 115, Taiwan}
\email{tonysun1981@gmail.com}

\author{Julie Tzu-Yueh Wang}
\address{Institute of Mathematics\\
 Academia Sinica\\
 6F, Astronomy-Mathematics Building\\
 No. 1, Sec. 4, Roosevelt Road \\
 Taipei 10617\\
 Taiwan}
\email{jwang@math.sinica.edu.tw}

\thanks{2020\ {\it Mathematics Subject Classification}: Primary 11J97; Secondary 14H05 and 11J87}
\thanks{The  first-named author was supported in part by National Natural Science Foundation of China (No. 12201643) and Natural Science Foundation of Hunan Province, China (No. 2023JJ40690).}
\thanks{The second-named author was supported in part by  NSERC grant RGPIN-2018-03770
and CRC tier-2 research stipend 950-231716 from the Government of Canada.}
\thanks{The   fourth-named author was supported in part by Taiwan's NSTC grant  110-2115-M-001-009-MY3.}

\begin{abstract}
We prove  Vojta's generalized abc conjecture for algebraic tori over function fields with  exceptional sets that can be determined effectively.  Additionally, we establish a version of the conjecture for toric varieties.  
As an application, we  investigate the Lang-Vojta Conjecture for  varieties of log general type that are ramified covers of $\mathbb G_m^n$ over function fields.   In particular, we consider the case  of $ \mathbb P^n\setminus D$, where $D$ is an algebraic curve over a function field in $\mathbb P^n$ with $n+1$ irreducible components and $\deg D\ge n+2$.  Our methods also apply to the complex situation, enabling us to find explicit exceptional sets for the corresponding case of Vojta's general abc conjecture (complex version) and the Green-Griffith-Lang conjecture.
\end{abstract}

\maketitle
\baselineskip=16truept

 
\section{Introduction }\label{sec:intro}
The primary objective of this article is to examine specific instances related to Vojta's general abc conjecture and Lang-Vojta conjecture in the context of function fields. We first present the conjectures.
We recall the notion of an ``admissible pair" $(X,V)$ introduced by  Levin in  \cite{Levin:GCD}, where $V$ is a nonsingular variety embedded in a nonsingular projective variety $X$ such that $D_0=X\setminus V$ is a normal crossings divisor.  Motivated by  \cite[Conjecture 1.13]{Levin:GCD}, we  reformulate Vojta's general abc conjecture   for   admissible pairs $(X,V)$.
\begin{conjecture} \label{Pairgeneralabc} 
Let $(X,V)$ be an admissible pair and let $D$ be a normal crossings divisor on $(X,V)$.  Let $A$ be a big divisor on $X$. 
Let $k$ be a number field over which $X$, $V$, $D$, and $A$ are defined and let $S$ be a finite set of places of $k$ containing all the archimedean places of $k$.  Then for any $\epsilon>0$, there is a proper Zariski-closed subset $Z$ of $X$, depending only on $X$, $D$, $ A$, and $\epsilon$ such that   for any set $R$ of $S$-integral points on $V$
\begin{align}\label{gabc}
N^{(1)}_{D,S}(x)\ge h_{{\mathbf K}_{(X,V)}+D}(x)-\epsilon h_{ A}(x)-O(1)
\end{align}
 holds for   all $x\in R\setminus Z$.
\end{conjecture}
Here,  $N_{D,S}^{(1)}$ is the truncated counting function with respect to $D$ and $S$,  $h_{ A}$ is the Weil height function associated with the divisor $A$, and
${\mathbf K}_{(X,V)}:={\mathbf K}_X+D_0$. By saying that $D$ is a normal crossings divisor on $(X,V)$, we mean that  $D$ is an effective divisor on $X$ such that $D+D_0$ is a normal crossings divisor on $X$.

We now recall the Lang-Vojta conjecture (\cite[Proposition 15.9]{Vojta}):
\begin{conjecture}\label{Lang-Vojta}
Let $k$ be a number field and $S$ be a finite set of places of $k$ containing all the archimedean places of $k$.
Let $  X$ be a smooth  projective variety over  $k$, and let $D$ be a   normal crossings divisor on $ X$ over $k$.  Let   ${\mathbf K}_X$   be a canonical divisor of $ X$.
If $V:=  X\setminus D$ is of log general type  (i.e.   ${\mathbf K}_X  +D$ is big), then no set of $(D,S)$-integral  points on $X(k)$ is Zariski dense.
\end{conjecture}

Let ${\bf k}$ be an algebraically closed field of characteristic zero,
$C$ be a smooth projective curve of genus $\mathfrak{g}$ defined over ${\bf k}$,
and $K:={\bf k}(C)$ be the function field of $C$.  Let $S$ be a finite set of points of $C({\bf k})$.  
Our primary result is the following function field version of Conjecture \ref{Pairgeneralabc} for an admissible pair, which consists of a toric variety   (see \cite[Definition 3.1.1]{cox})  and $\mathbb G_m^n$.
\begin{theorem}\label{toric}
Let $X$ be a nonsingular projective toric variety, and let $D_0:=X\setminus \mathbb G_m^n$.  Let $D$ be a normal crossings divisor on the admissible pair $(X,\mathbb G_m^n)$.  Let $A$ be a big divisor on $X$.
Let $X$, $D_0$, $D$ and $A$  be defined over $K$.  Let $\epsilon >0$.  Then there exists a proper closed subset $Z$ of $X$,   depending only on $X$, $D$, $ A$, and $\epsilon$  such that 
\begin{align}\label{gabctoric}
N_{D,S}^{(1)}(P)\ge h_D(P)-\epsilon h_A(P)-O(1)
\end{align}
for all $P\in \mathbb G_m^n(\mathcal O_S)\setminus Z$.
 \end{theorem}
\begin{remark*} In this case,   $D_0:=X\setminus \mathbb G_m^n$ is  a normal crossings divisor of $X$  by \cite[Theorem 3.1.19]{cox}  and \cite[Example 1.2.21]{cox}.   
Furthermore, ${\mathbf K}_{(X,V) }={\mathbf K}_X+D_0$ is trivial, as a canonical divisor ${\mathbf K}_X$ can naturally be taken  to be $-D_0$.(See \cite[Theorem 8.2.3]{cox}.)  Therefore, \eqref{gabctoric} implies \eqref{gabc}.
\end{remark*}

 The function field analogue of  Lang-Vojta conjecture has been studied by Corvaja-Zannier in the split case   for the complement of a conic and two lines in $\mathbb P^2$ in \cite{CZ2008}, 
for surfaces that are ramified covers of $\mathbb G_m^2$ in \cite{CZ2013}, and  by   Turchet \cite{Tur} and Capuano-Turchet \cite{CT}  in the corresponding non-split case. We  refer to  \cite{CT} and \cite{RTW} for  more discussions on related conjectures and results,   and \cite{Vol94} for a case of ramified covers of abelian varieties over function fields.  

We now present our results in the direction of Lang-Vojta conjecture for function fields.  We will start
 by considering the case where  $ X=\mathbb P^n$ and    $D$ is  a divisor of $\mathbb P^n$ over $K$ with $n+1$ irreducible components.   It is worth noting that  $\mathbb P^n\setminus D$ is of log general type if  $\deg D\ge n+2$.

\begin{theorem}\label{thmVconj}
Let $D=D_1+\cdots+D_{n+1}$, where $D_1,\hdots,D_{n+1}$ are hypersurfaces in $\mathbb P^n$ defined by irreducible homogeneous polynomials $ F_1,\hdots,F_{n+1}\in K[x_0,\hdots,x_n]$   with   $\sum_{i=1}^{n+1}\deg F_i\ge n+2$.  Suppose that there exists a point $\p\in C({\bf k}) $ such that the specializations $D_i(\p)$  of $D_i$, $1\le i\le n+1$,  intersect transversally.
Then there exists  
   a proper  Zariski closed set  $W\subset \PP^n$ over $K$ such that for any $(D,S)$-integral  subset $\Sigma\subseteq \PP^n(K)$ there is a constant $c=c(D,\Sigma)$  such that
\begin{equation}\label{htbounded}
h(\mathbf{x})\le c\max\{1,2\gen-2+|S|\} 
\end{equation}
 for any ${\mathbf x}\in \Sigma\setminus W$.
 Furthermore, if  $ F_i\in {\mathbf k}[x_0,\hdots,x_n]$ for $1\le i\le n+1$, then $W$ is defined over $\mathbf k$.
 \end{theorem}

The definition of a specialization  $R(\p)$  of  a hypersurface $R=[F=0]$ in $\mathbb P^n$  over $K$ will be given in Section \ref{Weilfunction}.   In particular, if  all coefficients of $F\in  K[x_0,\hdots,x_n]$ lie in $O_S^*$, then the specialization $R(\p)$ is the zero locus of  $F(\p)\in{\bf k}[x_0,\hdots,x_n]$ (i.e. evaluating the coefficients of $F$ at $\p$) in $\mathbb P^n({\bf k})$.
 
 \begin{remark*}  
The set $W$ can be determined explicitly and the degree of $W$ can  bounded  by an effectively computable constant by Theorem \ref{main_thm}.
\end{remark*}
 

\begin{theorem}\label{vojtaconj}
Let  $ X$ be a smooth  projective  variety over $K$  of dimension $n\ge 2$, $D$ be a  normal crossings divisor on $ X$, and $V= X\setminus D$. 
Suppose there is a morphism $\tilde \pi:    X\to\mathbb P^n$ defined over $K$   such that the restriction $\pi=\tilde \pi|_V : V\to \mathbb G_m^n$ is a finite morphism. 
Let $ Z\subset V$ be the ramification divisor of $\pi$, and let $R$ be the closure of $\pi(Z)$ in $\mathbb P^n$.  Additionally, assume  there exists a point $\p\in C({\bf k})$ such that the specialization $R(\p)\in \mathbb P^n({\bf k}) $ does not intersect the set of points  $\{(1,0,\hdots,0),\hdots,(0,\hdots,0,1)\} $. 
  Let $A$ be a big divisor on $ X$ over $K$.
If $V$ is of log-general type, then  there exists 
a proper Zariski closed subset $W$ of $ X$ over $K$  such that for any $(D,S)$-integral  subset  $\Sigma\subseteq  X(K)$ there exists a positive real  $c=c(\tilde \pi, \Sigma,A)$ satisfying
\begin{align}\label{htboundedA}
h_A({\mathbf x})\le c \max\{1,2\gen-2+|S|\}  
\end{align}
for any  ${\mathbf x}\in\Sigma\setminus W$.  
\end{theorem}
 \begin{remark*} 
The case $n=2$ has been proved by Corvaja and Zannier  in \cite{CZ2013}  for the split case, and  by Capuano and Turchet in  \cite{CT}  for the non-split case with the conclusion  that  $\deg \phi(C)$ is bounded in place of \eqref{htboundedA}.  Furthermore, if $X$, $D$,  $\tilde \pi$ and $A$ are all defined over $\mathbf k$, then $W$ is defined over $\mathbf k$.
\end{remark*}

\begin{remark*} Our methods are also applicable in the complex setting of Theorem \ref{thmVconj} and  Theorem \ref{vojtaconj}, allowing us to identify an exceptional set for the corresponding cases of the Green-Griffith-Lang conjecture established by Noguchi-Winkelman-Yamanoi in  \cite{noguchi2007degeneracy} (see Theorem \ref{GG_conj} and Theorem \ref{ramified} in Section \ref{ComplexRemark}). Notably, the existence of an exceptional set in \cite{noguchi2007degeneracy} had only been demonstrated for the case of $n=2$.
 \end{remark*}

The proofs of the aforementioned main theorems are based on a detailed estimation of multiple zeros of a polynomial in $n$ variables over $K$ when evaluated at $S$-unit arguments. The following is the technical theorem of bounding multiple zeros and a case of Vojta's generalized abc conjecture. We refer the reader to Section~\ref{sec:Preliminaries} for the meaning of the notations used in the theorem.


 \begin{theorem}\label{main_thm}
Let  $G\in K[ x_1,\hdots, x_n]$ be a non-constant    polynomial  with neither monomial  factors nor repeated factors. 
Then for any $\epsilon>0$,   there exist a positive real number $c_0$  and  
a proper Zariski closed subset $Z$ of $\mathbb{A}^n(K)$ such that for all  $(u_1,\hdots,u_n)\in ({\cal O}_{S}^*)^{n}\setminus Z$, we have either
\begin{enumerate}
\item the inequality 
\begin{align}\label{htubd} 
 \max_{1\le i\le n}\{ h(u_i)\}\le  c_0  \left(\tilde h(G) + \max\{1,2\gen-2+|S|\}\right) 
\end{align}
holds, 
  \item or the following two statements hold

 \begin{enumerate}
 \item[{\rm(a)}]  $N_{S}( G(u_1,\hdots,u_n) )-N_{S}^{(1)}( G(u_1,\hdots,u_n))\le \epsilon \max_{1\le i\le n}\{ h(u_i)\}$ if $G(u_1,\hdots,u_n)\ne 0$.
 \item[{\rm(b)}]  If $G(0,\hdots,0)\ne 0$ and $\deg_{X_i}G=\deg G=d$ for $1\le i\le n$, then 
 \begin{align*} 
 N_{S}^{(1)}( G(u_1,\hdots,u_n))\ge  \deg G \cdot (1-\epsilon)\cdot h(1,u_1,\hdots,u_n).
 \end{align*}
  \end{enumerate}
  \end{enumerate}
 Here, $c_0$ can be  effectively bounded from above in terms of $\epsilon$, $n$, and the degree of $G$.    Moreover, the exceptional set $Z$ can be expressed as the zero locus of a finite set
$\Sigma\subset K[x_1,\ldots,x_n]$ with the following properties: {\rm(Z1)} $\Sigma$ depends on $\epsilon$ and $G$ and can be determined explicitly,  {\rm(Z2)} $\vert \Sigma\vert$ and the degree of each polynomial in $\Sigma$ can be effectively bounded from above in terms of $\epsilon$, $n$, and the degree of $G$,  {\rm(Z3)}  if  $G\in \mathbf{k}[x_1,\ldots,x_n]$, then  we may take $\Sigma\subset \mathbf{k}[x_1,\ldots,x_n]$, and  {\rm(Z4)} if $n=2$, then $Z$ is a finite union of translates of proper algebraic  subgroups of $\mathbb G_m^2$.   
 \end{theorem}

 \begin{remark*}
 When $n\ge 3$, the exceptional set $Z$ may include Zariski closed subsets that are not   translates of  algebraic  subgroups of $\mathbb G_m^n$. The following are examples by taking $n=3$, $K={\bf k} (t)$, and $S=\{0,\infty\}$.
 \begin{enumerate}
 \item[{\rm(1)}] 
Let  $G= 1+2x_1-x_2+x_3$.  
 Let  $(u_1,u_2,u_3)=(t^{\ell}, at^{ 2\ell}, (1+a)t^{ 2\ell})$, where  where $\ell$ is any positive number  and $a\in \mathbf k\setminus\{0,-1\} $.
 Then $G(u_1,u_2,u_3)=(1+t^{\ell})^2$.  Then  $N_S(G(u_1,u_2,u_3))=2\ell$ and $N_S^{(1)}(G(u_1,u_2,u_3))=\ell$. Hence $$N_S(G(u_1,u_2,u_3))-N_S^{(1)}(G(u_1,u_2,u_3))=\ell=\frac12 \max\{h(u_1),h(u_2), h(u_3)\}.$$
 \item[{\rm(2)}] 
 Let $G= 1+x_1+x_2+x_3$.    Let  
  $(u_1,u_2,u_3)=(at^{\ell}, bt^{\ell}, ct^{\ell})$, where $\ell$ is any positive integer  and $(a,b,c)\in {\mathbf k^*}^3$ is an arbitrary tuple such that $a+b+c=0$.  Then  $G(u_1,u_2,u_3)=1$.  Hence $N_{S}^{(1)}( G(u_1,u_2,u_3))=0$.
   \end{enumerate}
\end{remark*}

%

To see how Theorem \ref{main_thm} is related to  Vojta's generalized abc conjecture and to illustrate the essential difficulties involved in the function field case, we state the complex version of the conjecture.  
	(See \cite[Conjecture 15.2]{Vojta}   and \cite[Conjecture 27.5]{Vojta}.)
	\begin{conjecture}\label{ConjABC}
		Let $X$ be a smooth complex projective variety,   D be a normal crossing divisor on $X$,  ${\mathbf K}_X$ be a canonical divisor on $X$, and  $A$ be an ample divisor on $X$. Then
		\begin{enumerate}
			\item[{\rm (a)}] If $f:\mathbb C\to X$ is an algebraically nondegenerate analytic map, then
			\begin{align}\label{truncate1_1}
				N_f^{(1)}(D,r)\ge_{\rm exc} T_{{\mathbf K}_X+D,f}(r)-{\rm o}(T_{ A,f}(r)).
			\end{align}
			\item[{\rm (b)}] For any $\epsilon>0$, there exists a proper Zariski-closed subset $Z$ of $X$, depending only on $X$, $D$, $ A$, and $\epsilon$ such that  for any non-constant analytic map $f:\mathbb C\to X$ whose image is not contained in $Z$, the following
			\begin{align}\label{truncate2}
				N_f^{(1)}(D,r)\ge_{\rm exc} T_{{\mathbf K}_X+D,f}(r)-\epsilon T_{ A,f}(r)
			\end{align}
			holds.
		\end{enumerate}
	\end{conjecture}
	Here, for each positive integer $n$, $N_f^{(n)}(D,r)$ is the $n$-truncated counting function with respect to $D$   given by
	\begin{align}
		N_f^{(n)}(D,r)=\sum_{0<|z|<r}\min\{{\rm ord}_z f^*D, n\}\log\frac{r}{|z|}+\min\{{\rm ord}_0 f^*D, n\}\log r,
	\end{align}
	where $T_{D,f}(r)$ is the (Nevanlinna) height function relative to the divisor $ D$ (referring to \cite[Section 12]{Vojta}), and the notation $\geq_{\rm exc}$  means that the estimate holds for all $r$ outside a set of finite Lebesgue measure.
	
Let $X$ be a smooth projective variety defined over the function field $K$.	
Let $\mathcal X$ be  a  proper model  of $X$, i.e. it is a  normal  variety $\mathcal X$  equipped with a proper flat morphism $\rho:\mathcal X\to C$, such that the generic fiber is isomorphic to $X$(see \cite[Section 16]{Vojta}).  Then the rational points in $X(K)$ correspond bijectively  to sections $i:C\to \mathcal X$.  When $X$ is defined over ${\bf k}$,  then the rational points in $X(K)$ correspond bijectively  to morphisms $i:C\to X$.  We will refer to the split case if $X$ is define over ${\bf k}$  and the non-split case otherwise.  Since every rational point in $X(K)$ is algebraically degenerate, the appropriate formulation of Vojta's general abc conjecture in the context of function fields can only be analogous to Conjecture \ref{ConjABC}(b).  Additionally, Conjecture \ref{ConjABC}(b) is analogous to the split case as $X$ and the divisors are over $\mathbb C$, whereas the non-split case should be viewed as instances of moving targets, i.e. the divisor $D$ is defined by meromorphic functions.  Consequently, addressing these scenarios within the framework of function fields is generally more challenging and subtle.
There are many results in this direction with a high truncation level in both complex  and  function field cases. However, the number of results available at level one is quite limited. Furthermore, the capacity to construct an explicit exceptional set remains highly constrained.   The following are some known results.
	First, Conjecture \ref{ConjABC} holds  for $\dim X=1$ (see \cite[Theorem 23.2]{Vojta}).  When $X$ is a semiabelian variety, Noguchi, Winkleman and Yamanoi in \cite{noguchi2008semiabelian} showed that the inequality \eqref{truncate2} holds if the map is algebraically nondegenerate.  The conjecture is much harder for the complex case with moving targets. 	The only existing results   in the case with slowly growth moving targets   are due to  Yamanoi in \cite{Yamanoi2004} for   $\dim X=1$, and in \cite {GSW22}, where  the inequality \eqref{truncate2} is derived for complex tori under the assumption that  the map is multiplicatively independent over a field of ``small functions" with respect to $f$. 

 In the function field situation (of characteristic zero),  a version of Vojta's general abc conjecture is proved by  Yamanoi for $\dim X=1$  in \cite{Yamanoi2004}  and   for the  case that $X$ is defined over the constant field  with maximal Albanese dimension in \cite{Yamanoi2015}.  To the best of our knowledge, there is no description of the exceptional sets in any of the aforementioned results, and the only related article with explicit description of the exceptional set  in this direction appears in  \cite{Garcia} with $ X=\mathbb P^2$ and $D$ consisting of  sufficiently many lines.

To see how Theorem \ref{main_thm} is related to Conjecture \ref{ConjABC}, we
let $\tilde G$ be the homogenization of $G$ in $K[x_0,\hdots,x_n]$, and let $D=[\tilde G=0]+[x_0=0]+ \cdots +[x_n=0]$ be a divisor of  $\mathbb P^n$. Then ${\mathbf K}_{\mathbb P^n}+D$ is linearly equivalent to $[\tilde G=0]$ and hence $h_{{\mathbf K}_{\mathbb P^n}+D}(x)=\deg G\cdot h(x)+O(1)$.  
If $\mathbf{u}=[u_0:\cdots:u_n]$ with $(u_0,\hdots,u_n)\in ({\cal O}_{S}^*)^{n+1}$, then we have $ N_{D,S}^{(1)} ( \mathbf{u} )=N_{S}^{(1)}( \tilde G(\mathbf{u})) +O(1)$.  Therefore, Theorem \ref{main_thm}    implies \eqref{truncate2} in the function field situation for  such $\mathbf{u}$, which can be identified as the point $(\frac{u_1}{u_0},\hdots, \frac{u_n}{u_0})$  in $\mathbb G_m^n({\cal O}_{S})$.

Our proof of Theorem \ref{main_thm} extends the work of Corvaja and Zannier in    \cite{CZ2008}, where they  derived the inequality stated in Theorem \ref{main_thm} (a) for $n=2$.  
The primary components of \cite{CZ2008} involve estimating the greatest common divisor of two multivariable polynomials evaluated at $S$-unit arguments over the field ${\bf k}$. Additionally, they introduce a method to bound the number of multiple zeros of such evaluations outside of $S$. These techniques are further refined and expanded upon for higher-dimensional cases in \cite{GSW}. Furthermore, in the same work, they establish a GCD theorem for multivariable polynomials over function fields, adopting the approach introduced by Levin in the paper \cite{Levin:GCD}. Nevertheless, this version of the GCD theorem comes with additional conditions that restrict its applicability to more general cases. In order to address these limitations and provide explicit descriptions of the exceptional sets required for the proof of Theorem \ref{main_thm}, we will derive specific specialization lemmas. These lemmas will be crucial for conducting the proof by induction.

 In Section \ref{sec:Preliminaries}, we will review fundamental definitions of heights and integral points, state a modified version of a lemma of Brownawell and Masser on $S$-unit equations over function fields and establish essential properties for polynomials related to $n$-tuples of units. The technical lemmas of specialization will be presented in Section \ref{mainlemmas}. We will prove   Theorem \ref{main_thm} in Section \ref{Theorem1general} and Theorem \ref{toric} in  Section \ref{Theorem2general}. In Section \ref{LangVojtaF}, we will demonstrate the proofs for Theorem \ref{thmVconj} and Theorem \ref{vojtaconj}.  Finally, in Section \ref{ComplexRemark}, we will briefly discuss, without presenting a full proof, the application of the reduction method in proving Theorem \ref{main_thm} in the complex case. This will enable us to derive explicit exceptional sets for both Vojta's general abc conjecture (complex version) and the (strong) Green-Griffith-Lang conjecture corresponding to Theorem  \ref{toric}- \ref{vojtaconj}.  
\section{Preliminaries}\label{sec:Preliminaries}
\subsection{Notation and heights}\label{Notation and heights}
Let ${\bf k}$ be an algebraically closed field of characteristic zero,
$C$ be a smooth projective curve  of genus $\mathfrak{g}$ defined over ${\bf k}$,
and $K:={\bf k}(C)$ be the function field of $C$.
 At each point $\p\in C(\mathbf{k})$,
we may choose a uniformizer $t_{\p}$ and define the normalized order
function $v_{\p}:=\ord_{\p}:\K\to\ZZ\cup\{+\infty\}$. 
We  denote by 
$$
\chi_S(C):= 2\mathfrak{g}-2+|S|, \quad \text{and} \quad \chi_S^+(C):=\max\{0, \chi_S(C)\}.
$$
For $f\in\K^{*}$ and $\mathbf{p}\in C(\mathbf{k})$
we let 
\[
v_{\p}^{0}(f):=\max\{0,v_{\p}(f)\},\quad \text{and}\quad v_{\p}^{\infty}(f):=-\min\{0,v_{\p}(f)\}
\]
i.e. its order of zero and poles at $\p$ respectively.
 The height
of $f$ is defined by 
\[
h(f):=\sum_{\p\in C(\mathbf{k})}v_{\p}^{\infty}(f).
\]
 For
any ${\bf f}:=[f_{0}:\cdots:f_{n}]\in\PP^{n}(K)$ with $n\ge1$ and
$f_{0},...,f_{n}\in\K$, we define 
$$v_{\p}(\mathbf{f}):=\min\{v_{\p}(f_{0}),...,v_{\p}(f_{n})\}$$
and 
\[
h({\bf f})=h(f_{0},...,f_{n}):=\sum_{\p\in C(\mathbf{k})} -v_{\p}(\mathbf{f}).
\]
For a finite subset $S$ of $C(\mathbf{k})$, we define the ring of $S$-integers and group of $S$-units in $K$ by
$${\cal O}_{S}:=\{f\in\K\,|\,v_{\p}(f)\ge0\text{ for all }\p\notin S\},\quad\text{and}\quad {\cal O}_{S}^{*}:=\{f\in\K\,|\,v_{\p}(f)=0\text{ for all }\p\notin S\}.$$
For  $f\in K^{*}$ and  a positive integer $m$, we let 
\[
N_{S}({f})=\sum_{\mathbf{p}\in C({\bf k})\setminus S} v_{\p}^{0}(f) \quad \text{and}\quad {N}^{(m)}_{S}({f})={\displaystyle \sum_{\mathbf{p}\in C({\bf k})\setminus S}\min\{m,{v}_{\mathbf{p}}^{0}(f)}\}
\]
be the number of the zero of $f$ outside of $S$, counting multiplicities and counting multiplicities up to $m$ respectively.
 For any $f,g\in\K,$ we let 
\begin{align*}
N_{S,{\rm gcd}}(f,g):=\sum_{\p\in C({\bf k})\setminus S}\min\{v_{\p}^{0}(f),v_{\p}^{0}(g)\}\quad \text{and}\quad h_{{\rm gcd}}(f,g):=\sum_{\p\in C({\bf k})}\min\{v_{\p}^{0}(f),v_{\p}^{0}(g)\}.
\end{align*}

Let $\mathbf{x}:=(x_{1},\ldots,x_{n})$ be a tuple of $n$ variables,
and $F=\sum_{{\bf i}\in I_{F}}a_{{\bf i}}{\bf x}^{{\bf i}}\in K[x_{1},\ldots,x_{n}]$
be a nonzero polynomial, where $I_{F}$ is the set
of those indices ${\bf i}=(i_{1},\hdots,i_{n})$ with $a_{{\bf i}}\ne0$;
and we put ${\bf x}^{{\bf i}}:=x_{1}^{i_{1}}\cdots x_{n}^{i_{n}}$.
We define the height $h(F)$ and the relevant height $\tilde{h}(F)$
as follows. Put 
\begin{align*} 
v_{\p}(F):=\min_{{\bf i}\in I_{F}}\{v_{\p}(a_{{\bf i}})\}\qquad\text{for }  \p\in C({\bf k}),
\end{align*}
and define 
\begin{align*}
h(F):=\sum_{\p\in C({\bf k})}-v_{\p}(F), \qquad\text{and }\quad \tilde{h}(F):=\sum_{  \p\in C({\bf k})}-\min\{0,v_{\p}(F)\}.
\end{align*}

Notice that Gauss's lemma can be stated as 
\begin{align*}
v_{\p}(FG)=v_{\p}(F)+v_{\p}(G), 
\end{align*}
where $F$ and $G$ are in $K[x_{1},\hdots,x_{n}]$ and $\p\in C({\bf k})$.
Consequently, we have that 
\begin{align}\label{hp}
h(FG)=h(F)+h(G),\qquad\text{and }\quad \tilde h(FG)\le \tilde h(F)+\tilde h(G).
\end{align} 
 
We will use the following slightly modified result of Brownawell-Masser \cite{BM}.
\begin{theorem}\label{BrMa}  
If $f_{1},\hdots,f_{n}\in K^* $ and $f_{1}+\cdots+f_{n}=1$,
then either some proper subsum of $f_{1}+\cdots+f_{n}$ vanishes or
\[
 h(1,f_{1},\hdots,f_n)\le \sum_{i=1}^n \big(N_{S}^{(n)}(f_i)+N_{S}^{(n)}(f_i^{-1})\big)+\frac{n(n-1)}{2}\max\{0,2\gen-2+|S|\},
\]
where $S$ is any finite subset of $C$.
\end{theorem} 
The proof can be deduced from the arguments in \cite{BM}, or  from adapting the proof in the complex version  \cite[Theorem 2.1]{RuWang2003}.
\begin{corollary}\label{ProximityAffine} 
Let $F=\sum_{\mathbf i\in I_F}a_{\mathbf i}{\mathbf x}^{\mathbf i}\in K[ x_1,\hdots, x_n]$ be a non-constant polynomial.  Assume that $F(0,\dots,0)\ne 0$.  Let $Z$ be the Zariski closed subset that is the union of hypersurfaces of $\mathbb A^n$ of the form  $\sum_{\mathbf i\in I}\alpha_{\mathbf i}{\mathbf x}^{\mathbf i}=0$ where $I$ is a non-empty subset of $I_F$.
 Let $S$ be a finite subset of $C({\mathbf k})$.  Then
for all $(u_1,\hdots,u_n)\in (\mathcal O_S^*)^n\setminus Z$, we have 
 \begin{align*} 
 \sum_{\p\in S} v_{\p}^0(F(u_1,\hdots,u_n)) \le  
\tilde c_1 \chi_S^+(C) +\tilde c_2h(F),
  \end{align*}
  where $d=\deg F$, $\tilde c_1=\frac12 \binom{n+d}{n}  (\binom{n+d}{n}+1)$ and $\tilde c_2=2 (\binom{n+d}{n}-1)\tilde c_1$.
\end{corollary}
 \begin{proof}
Since $F(0,\hdots,0)\ne 0$ and $h(F)=h(\lambda F)$ for $\lambda\in K^*$, we may assume that $F(x_1,\hdots,x_n)=1+\sum_{\mathbf i\in I'_F}a_{\mathbf i} \mathbf{x}^{\mathbf i}$, where $\mathbf{i}\in\mathbb{Z}_{\ge0}^{n}$,
 $|\mathbf{i}|\le d$, and $a_{\mathbf i}\ne 0$.  Then we have
 \begin{align}\label{Fequ}
 1=F( {\mathbf u})-\sum_{\mathbf i\in  I'_F}a_{\mathbf i}  {\mathbf u}^{\mathbf i},
 \end{align}
where $ {\mathbf u}=(u_1,\hdots,u_n)\in (\mathcal O_S^*)^n\setminus Z$.
Therefore,  no subsum of $\sum_{\mathbf i}a_{\mathbf i} {\mathbf u}^{\mathbf i}$ vanishes.  
Let
 $S'=\{\p\in C\setminus S\,|\, v_{\p}(a_{\mathbf i})\ne 0 \text{ for some } {\mathbf i\in I'_F} \}\cup S.$  Then   $F( {\mathbf u})\in\mathcal O_{S'}$ and
 $$
 |S'|\le |S|+2c_1 h(F),
 $$
 where $c_1=  \binom{n+d}{n}-1$.
We now apply Theorem \ref{BrMa} to \eqref{Fequ}  with $S'$ to get
 \begin{align}\label{htbdd}
 h(F( {\mathbf u})) \le N_{0,S}(F( {\mathbf u})) 
  +\tilde c_1  \chi_S^+(C) +2c_1\tilde c_1h(F)
 \end{align}
where $\tilde c_1=\frac12 (c_1+1)(c_1+2)$.
Since 
$h(F( {\mathbf u}))=\sum_{\p\in C} v_{\p}^0(F( {\mathbf u})),$
  the inequality \eqref{htbdd} implies
 \begin{align*} 
 \sum_{\p\in S} v_{\p}^0(F( {\mathbf u})) \le  
\tilde c_1 \chi_S^+(C) +2c_1\tilde c_1h(F).
 \end{align*}
\end{proof}

\subsection{Weil functions associated to divisors and integral points}\label{Weilfunction}
 We now recall some facts from \cite[Chapter 10]{LangDG} and \cite[Section B.8]{HS} about local Weil functions associated to divisors.   
We continue to let $K$ be the function field of a smooth projective curve $C$ over   an algebraically closed field ${\bf k}$  of characteristic zero.  Denote by  $M_K:=\{v=v_{\p}:\p\in C(\mathbf{k})\}$ the set of  valuations on $K$.
We recall that an \emph{$M_K$-constant} is a family $\{\gamma_v\}_{v\in M_K}$, where each $\gamma_v$ is a real number with all but finitely many being zero.  
Given two families $\{\lambda_{1v}\}$ and $\{\lambda_{2v}\}$ of functions parametrized by $M_K$, we say $\lambda_{1v} \le \lambda_{2v}$ holds up to an $M_K$-constant if there exists an $M_K$-constant $\{\gamma_v\}$ such that the function $\lambda_{2v} - \lambda_{1v}$ has values at least $\gamma_v$ everywhere.  We say $\lambda_{1v} = \lambda_{2v}$ up to an $M_K$-constant if $\lambda_{1v} \le \lambda_{2v}$ and $\lambda_{2v} \le \lambda_{1v}$ up to $M_K$-constants. Let $X$ be a projective variety over $K$.  
We say that a subset $Y$ of $X(K)\times M_K$ is {\it affine $M_K$-bounded} if there is an affine open subset $X_0\subset X$   over $K$ with   a system of affine coordinates $x_1,\hdots,x_n$ and an $M_K$-constant $\{\gamma_v\}_{v\in M_K}$ such that $Y\subset X_0 (K) \times M_K$ and  
$$
\min_{1\le i\le n}v(x_i(P))\ge  \gamma_v, \qquad\text{for all  } (P,v)\in Y;
$$
and we say that the set $Y$ is {\it  $M_K$-bounded}  if it is a finite union of affine $M_K$-bounded sets.

The classical theory of heights associates to every Cartier divisor $D$ on $X$ a \emph{height function} $h_D:X(K)\to \mathbb R$ and a \emph{local Weil function}  (or \emph{local height function})
$\lambda_{D,\p}: X(K)\setminus {\rm Supp}(D)\to \mathbb R$
for each $v_{\p}\in M_K$, well-defined up to  an $M_K$-constant, such that
\begin{align}\label{FMT}
\sum_{\p\in C(\mathbf{k})}\lambda_{D,\p}(P)=h_D(P)+O(1)
\end{align}
for all $P\in X(K)\setminus {\rm Supp}(D)$. 
In particular, for   a nonzero rational function $f$ on $X$ with $D:={\rm div}(f)$, the difference 
$$
\lambda_{D,v}(P)-v(f(P))
$$ 
is an $M_K$-bounded function on every $M_K$-bounded subset of $X(K)\setminus {\rm Supp}(D)$, i.e., on every such subset, this family $\{\lambda_{D,v}(P)-v(f(P))\}_{v\in M_K}$ of real functions is bounded in the  usual absolute value by an $M_K$-constant. (See  \cite[Theorem B.8.1 (a)]{HS}.)

 For a finite subset $S$ of $C(\mathbf{k})$, we will denote by 
 $$
 m_{D,S}(x):=\sum_{\p\in  S}\lambda_{D,\p}(x), 
$$
$$
N_{D,S}(x):=\sum_{\p\in C(\mathbf{k})\setminus S}\lambda_{D,\p}(x),\quad\text{and}\quad 
N^{(m)}_{D,S}(x):=\sum_{\p\in C(\mathbf{k})\setminus S}\min\{m,\lambda_{D,\p}(x) \}
$$ 
where  $m$ is a positive integer.

Let $F\in K[x_0,\hdots,x_n]$ be a non-constant homogeneous polynomial and $D$ be   the divisor associated with $F$ in $\mathbb P^n$.   We  always take the   following local Weil functions:
\begin{equation}
\lambda_{D,\p}(\mathbf x):= v_{\p}(F(\mathbf x))- v_{\p}(F)-\deg F\cdot v_{\p}(\mathbf x)\ge 0 \label{Weil_ex}
\end{equation}
for $\p\in C ({\bf k}) $ and $\mathbf x\in\mathbb P^n(K)\setminus {\rm Supp}(D) $.
We note that this definition is independent of the choice of the projective coordinates of $\mathbf x$ and is invariant by $K^*$-scalar multiplications on $F$.     
Furthermore, the specialization $D(\p)$ of $D$ at a point $\p\in C$ is simply the zero locus of the $ F_{\p}(\p)$ in $\mathbb P^n({\bf k})$, where $F_{\p}=t_{\p}^{-v_{\p}(F)}F$ and $F_{\p}(\p)\in{\bf k}[x_0,\hdots,x_n]$ is obtained from  $F_{\p}$ by evaluating its coefficients at $\p$. 

Let $D_i$, $1\le i\le q$,  be  hypersurfaces  in $\mathbb P^n$ defined over a field $\mathcal K$.  We say that $D_1,\hdots,D_q$ are in general position if 
the intersection  in $\mathbb P^n(\overline {\mathcal K})$ of any distinct $n+1$ hypersurfaces among $D_i$, $1\le i\le q$, is empty.
The following proposition shows that the hypersurfaces over $K$ are in general position   if one of their specializations are in general position.  It can be proved easily by the theorem of resultant.(See  \cite[Chapter IX]{langalgebra} or \cite[Chapter XI]{waerden1967}. )

\begin{proposition}\label{generalposition}
Let $D_i$, $1\le i\le q$, be hypersurfaces  in $\mathbb P^n$ defined over $K$.  Suppose that there exists a point $\p\in C({\bf k})$ such that the specializations $D_i(\p)$ of $D_i$, $1\le i\le q$, are in general position.  Then $D_1,\hdots,D_q$ are in general position.\end{proposition}


Finally, we recall the following definition of $(D,S)$-integral   subsets of $X(K)$  for a projective variety  $X$ over $K$ and an effective divisor $D$ on $X$ over $K$  analogous to the number field setting.
\begin{definition}\label{int_set}
 Let $X$ be a projective variety over $K$, and let $D$ be an effective Cartier divisor on $X$ over $K$.   A $(D,S)$-integral  subset  $\Sigma\subseteq X(K)$ is a subset containing no point in ${\rm Supp}(D)$ such that there is a Weil function $\{\lambda_{D,\p}\}_{\p\in C({\bf k})}$ for $D$ and an $M_K$-constant $\{\gamma_\p\}$ satisfying that $\lambda_{D,\p}(P)\le \gamma_\p$ for all $\p\notin S$ and all $P\in \Sigma$. 
\end{definition}
 Any  Weil function for an effective Cartier divisor is bounded from below by an $M_K$-constant. Any two Weil functions for the same Cartier divisor differ by a function bounded by an $M_K$-constant. Thus, for every Weil function $\{\lambda_{D,\p}\}_{\p\in C({\bf k})}$ for an effective Cartier divisor $D$ on  a projective variety $X$  over $K$ and   every
$(D,S)$-integral subset of $\Sigma\subset X(K)$, there is {\em always} a finite subset $S'\subset C({\bf k})$ containing  $S$  such that $\lambda_{D,\p}(P)=0$ for each $P\in \Sigma$ and each $\p\notin S'$.

We need the following result, which follows from the basic properties of Weil functions. 
\begin{proposition}[{\cite[Proposition 13.5]{Vojta}}]\label{under_mor}
Let $\phi:X_1\to X_2$ be a $K$-morphism between projective varieties over $K$, and let $D_1, D_2$ be respectively effective Cartier divisors on $X_1,X_2$ over $K$. Assume that ${\rm Supp}(D_1)\supset { \rm Supp}(\phi^*(D_2))$. Then the image of a $(D_1,S)$-integral  subset of  $ X_1(K)$ under $\phi$ is a $(D_2,S)$-integral  subset of  $ X_2(K)$. 
\end{proposition}

\subsection{Polynomials associated with $n$-tuples of $S$-units}
We now recall the definitions of global and local derivations on $\K$.
Let $t\in\K\setminus {\bf k}$, which will be fixed later. The mapping
${\displaystyle {g\to\frac{dg}{dt}}}$ on ${\bf k}(t)$,  the  formal differentiation
on ${\bf k}(t)$ with respect to $t$, extends uniquely to a global
derivation on $\K$ as $\K$ is a finite separable extension of ${\bf k}(t)$.
Furthermore, since an element in $\K$ can be written as a Laurent
series in $t_{\p}$, the local derivative of $\eta\in\K$
with respect to $t_{\p}$, denoted by ${\displaystyle {d_{\p}\eta:=\frac{d\eta}{dt_{\p}}}}$,
is given by the formal differentiation on ${\bf k}((t_{\p}))$ with
respect to $t_{\p}$.   The chain rule says 
\begin{align}
\frac{d\eta}{dt}=d_{\p}\eta\cdot(d_{\p}t)^{-1}.\label{chain rule}
\end{align}

The following is a consequence  of the Riemann-Roch Theorem.
We refer to \cite[Corollary 7]{Buchi2013} for a proof. 
\begin{proposition}\label{functiont}
For a point $\q\in C({\bf k})$, we can find some $t\in K\setminus{\bf k}$ satisfying the following conditions:
$t$ has exactly one pole at $\q$; $h(t)\le\gen+1$; and ${\displaystyle {\sum_{\p\in C({\bf k})}v_{\p}^{0}(d_{\p}t)\le3\gen}}$.
 \end{proposition}

From now on, we will fix a $t$ satisfying the
conditions in Proposition \ref{functiont} and use the notation $\eta':=\frac{d\eta}{dt}$
for $\eta\in K$.

Next, we recall some definitions and results from \cite{GSW}.
For the convenience of discussions, we will use the following convention.
Let $\mathbf{i}=(i_{1},\ldots,i_{n})\in\mathbb{Z}^{n}$ and $\mathbf{u}=(u_{1},\ldots,u_{n})\in(K^{*})^{n}$.
We denote by $\mathbf{x}:=(x_{1},\ldots,x_{n})$, $\mathbf{x^{i}}:=x_{1}^{i_{1}}\cdots x_{n}^{i_{n}}$,
$\mathbf{u^{i}}:=u_{1}^{i_1}\cdots u_{n}^{i_{n}}\in K^{*}$ and $|\mathbf{i}|:=\sum_{j=1}^{n}|i_{j}|$.
For a polynomial $F(\mathbf{x})=\sum_{\mathbf{i}}a_{\mathbf{i}}\mathbf{x}^{\mathbf{i}}\in K[x_{1},\dots,x_{n}]$,
we denote by $I_{F}$ the set of exponents ${\bf i}$ such that $a_{\mathbf{i}}\ne0$
in the expression of $F$, and define 
\begin{align}
D_{\mathbf{u}}(F)(\mathbf{x}):=\sum_{\mathbf{i}\in I_{F}}\frac{(a_{\mathbf{i}}\mathbf{u}^{\mathbf{i}})'}{\mathbf{u}^{\mathbf{i}}}\mathbf{x}^{\mathbf{i}}.\label{DuF}
\end{align}
Clearly, we have 
\begin{align}
F(\mathbf{u})'=D_{\mathbf{u}}(F)(\mathbf{u}),\label{value}
\end{align}
and the following product rule: 
\begin{align}
D_{\mathbf{u}}(FG)=D_{\mathbf{u}}(F)G+FD_{\mathbf{u}}(G)\label{product}
\end{align}
for each $F,G\in K[x_{1},\dots,x_{n}]$.

We will use the following proposition and lemma from \cite{GSW}.
\begin{proposition}[{\cite[Proposition 17]{GSW}}]
\label{heightDu} Let $F$ be
a non-constant  polynomial in $K[x_{1},\dots,x_{n}]$ and $\mathbf{u}=(u_{1},\ldots,u_{n})\in(O_{S}^{*})^{n}$.
Then there exist $c_{1},c_{2}$ depending only on $\deg F$
such that 
\[
\tilde{h}(D_{\mathbf{u}}(F))\le c_{1}\tilde{h}(F)+c_{2}\max\{1,\chi_S(C)\}.
\]
\end{proposition}  
\begin{lemma} [{\cite[Lemma 18]{GSW}}]\label{lem:coprime-irr}
For any irreducible $F(\mathbf{x})=\sum_{\mathbf{i}\in I_{F}}a_{\mathbf{i}}\mathbf{x}^{\mathbf{i}}\in K[x_{1},\dots,x_{n}]$
and $\mathbf{u}\in(K^{*})^{n}$, the two polynomials $F$ and $D_{\mathbf{u}}(F)$
are not coprime if and only if $\frac{a_{\mathbf{i}}\mathbf{u}^{\mathbf{i}}}{a_{\mathbf{j}}\mathbf{u}^{\mathbf{j}}}\in{\bf k}^{*}$
whenever $\mathbf{i},\mathbf{j}\in I_{F}$. 
\end{lemma}
For our purpose, we adapt Lemma 20 in \cite{GSW} into the following formulation.
 \begin{lemma} \label{lem:coprime-gen} Let $F=\prod_{i=1}^{q}P_{i}\in K[x_{1},\hdots,x_{n}]$,
where $P_{i}$, $1\le i\le q$, are distinct irreducibles   in
$K[x_{1},\hdots,x_{n}]$. Let $\mathbf{u}\in(\mathcal O_S^{*})^{n}$. 
Then  we may factor $F=A \cdot B\in K[x_{1},\hdots,x_{n}]$ such that $B(\mathbf{u})\in \mathcal O_S^{*}\cup\{0\}$, $A$ is constant  or  $A$ and $D_{\mathbf{u}}(F)$  
are coprime in $K[x_{1},\hdots,x_{n}]$,  and $\tilde h(A)\le 2\tilde h(F)$. 
\end{lemma}
\begin{proof}
By the product  rule \eqref{product}, we have
\begin{equation}\label{diffrentialexpansion}
    D_{\mathbf{u}}(F) =D_{\mathbf{u}}(P_1) P_2\cdots P_q+P_1D_{\mathbf{u}}(P_2) P_3\cdots P_q+\cdots +P_1\cdots P_{q-1}D_{\mathbf{u}}(P_q).
\end{equation}
Then clearly $F$ and $D_{\mathbf{u}}(F)$ are coprime if 
 each pair $P_i$ and $D_{\mathbf{u}}(P_i) $, $1\le i\le q$, is coprime.   Therefore, we simply take $A=F$ and $B=1$.
Otherwise, we 
let $B$ be the product of those $P_i$ with the property that  $P_i$ is not coprime with $D_{\mathbf{u}}(P_i)$, and put $A=\frac{F}{B}$.   With respect  to a fixed monomial order on $K[x_{1},\hdots,x_{n}]$, we denote by ${\rm LC}(P)$, {\em the leading coefficient of a polynomial $P\in K[x_{1},\hdots,x_{n}]$},  the coefficient attached to the largest monomial appearing in $P$. Having the fact that ${\rm LC}(F)={\rm LC}(A){\rm LC}(B)$, we may rearrange so that ${\rm LC}(P_i)=1$ for each  $P_i$ dividing $B$ and ${\rm LC}(F)={\rm LC}(A)$. Denote by  $\a={\rm LC}(F)={\rm LC}(A)\in K^*$, $F_1=\frac 1{\alpha}F$ and $A_1=\frac 1{\alpha}A$, we see that $F_1=A_1 B$,  $ h(\alpha)\le \tilde h(F)$ and 
$$
\tilde h(A_1)=h(A_1)\le h(F_1)=h(F)\le  \tilde h(F).
$$
It is easy to check that
\begin{equation}\label{height_est}
\tilde h(A)\le h(\alpha)+\tilde h(A_1)\le 2\tilde h(F).
\end{equation}
 If $A$ is not constant, then by construction, each of those $P_i$ dividing $A$ is coprime with $D_{\mathbf{u}}(P_i)$, and thus does not divide $D_{\mathbf{u}}(F)$ by \eqref{diffrentialexpansion}.  Hence, $A$ is coprime to $D_{\mathbf{u}}(F)$ as desired.

It is left to show that $B(\mathbf{u})\in \mathcal O_S^{*}\cup\{0\}$.
 For each  $P_i$ dividing $B$, we have that   $P_i=\sum_{\mathbf{i}\in I_{P_i}}a_{\mathbf{i}}\mathbf{x^i}$ is not coprime with $D_{\mathbf{u}}(P_i)$ and  that $a_{\mathbf{j}}=1$ for some $\mathbf{j}\in I_{P_i}$.
By Lemma \ref{lem:coprime-irr}, we have $\frac{a_{\mathbf{i}}\mathbf{u^i}}{\mathbf{u^j}}\in {\bf k}^*$ whenever $\mathbf{i}\in I_{P_i}$.  
Thus we have $P_i(\mathbf{u})=c\mathbf{u^j}$ for some  $c\in {\bf k}$ and hence   $B(\mathbf{u})\in \mathcal O_S^{*}\cup\{0\}$ as desired. 
\end{proof}

\section{Some specialization Lemmas}\label{mainlemmas} 
\begin{lemma}\label{main_lemma}
Let $n\ge 2$ and let $(m_1,\hdots,m_n)$ be  a non-zero vector in $\mathbb{Z}^n$
with $\gcd(m_1,\ldots,m_n)=1$. Then 
there exist $\bfv_i=(v_{i,1},\ldots,v_{i,n})\in \mathbb{Z}^{n}$ for $1\leq i\leq n-1$ such that
$$\vert v_{i,j}\vert \leq \max\{\vert m_j\vert,1\}\quad \text{for  $1\leq j\leq n$}$$
and  
$(m_1,\hdots,m_n)$ together with the $\bfv_i$'s form a basis of  $\mathbb Z^n$.
\end{lemma}
\begin{proof}
This is given in \cite[pp.~175--177]{MTW08}.
\end{proof}

For the rest of this section, let $k$ be a field and let $q$ and $r$ be positive integers. We write $\bft$ and $\bfx$ to denote the tuples of indeterminates $(t_1,\ldots,t_q)$ and $(x_1,\ldots,x_r)$ respectively. For $\bfi=(i_1,\ldots,i_r)\in \mathbb{Z}^r$, denote $\bfx^{\bfi}=x_1^{i_1}\cdots x_r^{i_r}$. The notation $\bft^{\bfj}$ for $\bfj\in\mathbb{Z}^{q}$ is defined similarly. For $\Sigma\subseteq k[\bft]$, let $\mathcal{Z}(\Sigma)=\{\lambda\in k^q:\ f(\lambda)=0\ \text{for every $f\in \Sigma$.}\}$. If $k$ is infinite and $f\in k[\bft]\setminus\{0\}$ then we can easily prove
 $\mathcal{Z}(f)\subsetneq k^q$  by induction on $q$.

\begin{lemma}\label{lem:DPf}
Let $f(\bft,\bfx)\in k[\bft,\bfx]$ and let  $P(\bfx)\in k[\bfx]$ be non-zero. 
Put 
$$D_{P,f}=\{\lambda\in k^q:\ P(\bfx)\mid f(\lambda,\bfx)\},$$
then we have:
\begin{itemize}
\item [{\rm (a)}] Either $D_{P,f}=k^q$, 
\item [{\rm (b)}] or there exists an effectively computable finite non-empty set
$\Sigma\subset k[\bft]\setminus\{0\}$ such that $D_{P,f}=\mathcal{Z}(\Sigma)$. The cardinality of $\Sigma$ and the degree of each polynomial in $\Sigma$ can be bounded effectively in terms of $q$, $r$, and the degrees of $f$ and $P$.
\end{itemize}
Moreover, when $k$ is infinite then (a) happens if and only if $P(\bfx)\mid f(\bft,\bfx)$ in $k[\bft,\bfx]$. 
\end{lemma}
\begin{proof}
We fix a monomial ordering on the monials in $t_1,\ldots,t_q,x_1,\ldots,x_r$. Then we may take $\{P\}$ as a Gr\"obner basis of the principal ideal $(P)$ in $k[\bft,\bfx]$ and perform the general polynomial  division algorithm as given in \cite[p.~320]{DF03_AA}. Since $P\in k[\bfx]$ and $f\in k[\bft,\bfx]$, this algorithm yields
$$f=QP+R$$
with $Q,R\in k[\bft,\bfx]$ and $R$ is either $0$ or every monomial in $R$ is not divisible
by the leading monomial in $P$.

If $R=0$, we have $D_{P,f}=k^q$. Otherwise, write $R(\bft,\bfx)=\sum_{\bfi\in I}B_{\bfi}(\bft)\bfx^{\bfi}$ with
$I\neq \emptyset$ and $B_{\bfi}\neq 0$ for every $\bfi\in I$.
We have that $P\mid f(\lambda,\bfx)$ if and only if $R(\lambda,\bfx)=0$ \cite[pp.~321--322]{DF03_AA}. The latter is equivalent to $B_{\bfi}(\lambda)=0$ for every $\bfi\in I$. By following steps in the division algorithm in \cite[p.~320]{DF03_AA}, we can bound $\vert I\vert$ and the degree of each $B_{\bfi}$ explicitly in terms of $q$, $r$, and the degrees of $f$ and $P$.

For the last assertion, we have that the ``if'' part is immediate. The ``only if'' part follows from the above arguments: if $R\neq 0$ then the set $\mathcal{Z}(\{B_{\bfi}:\ \bfi\in I\})$ is strictly contained in $k^q$ since $k$ is infinite.
\end{proof}

\begin{lemma}\label{special} 
Assume that $k$ is infinite.
 Let $f(\bft,\bfx)\in k[\bft,\bfx]$ be a polynomial with no monomial factor and no repeated irreducible factor in $ k[\bft,\bfx]$.  Then there exists an effectively computable non-empty finite set
 $\Sigma\subset k[\bft]\setminus\{0\}$
 such that
 for every $\lambda\in k^q\setminus\mathcal{Z}(\Sigma)$, 
 the polynomial $f(\lambda,\bfx)$ has no monomial or repeated irreducible factor.
 Moreover, the cardinality of $\Sigma$ and the degree of each polynomial in $\Sigma$ can be bounded effectively in terms of $q$, $r$, and the degree of $f$.  
\end{lemma}
\begin{proof}
The case $f\in k[\bft]$ is trivial: we simply take $\Sigma=\{f\}$. 
We assume $f\notin k[\bft]$ from now on. We consider the sets:
$$\mathcal{S}_1=\{\lambda\in k^q:\ f(\lambda,\bfx)\ \text{has a monomial factor}\},$$
$$\mathcal{S}_2=\{\lambda\in k^q:\ f(\lambda,\bfx)\ \text{has a repeated irreducible factor}\}.$$

We have $\displaystyle\mathcal{S}_1=\bigcup_{P} D_{P,f}$ where $P$ ranges over the finitely many monomials in $x_1,\ldots,x_r$ whose degree is at most the degree of $f$. Then we apply Lemma~\ref{lem:DPf} to get $\mathcal{S}_1=\mathcal{Z}(\Sigma_1)$ for an effectively computable $\Sigma_1$ and the cardinality $\Sigma_1$ and
the degree of its elements can be bounded as in the conclusion of the lemma. It remains to treat $\mathcal{S}_2$.

Let $I=\{1\leq i\leq r:\ \deg_{x_i}(f)>0\}$ which is non-empty since $f\notin k[\bft]$.  For $i\in I$, let
$f_i=\displaystyle\frac{\partial f}{\partial x_i}$, let $k[\bft,\bfx]_i$ be the polynomial ring in the indeterminates $t_1,\ldots,t_q,x_1,\ldots,x_{i-1},x_{i+1},\ldots,x_r$ and let 
$$R_i=\Res_{x_i}(f,f_i)\in k[\bft,\bfx]_i$$
be the resultant of $f$ and $f_i$ regarded as polynomials in $x_i$ with coefficients in 
$k[\bft,\bfx]_i$. Since $f$ has no repeated factor, we have that $R_i\neq 0$. Regard $R_i$ as a polynomial over $k[\bft]$ and let $\Sigma_i'\subset k[\bft]\setminus\{0\}$ be the non-zero coefficients of $R_i$.
Write:
$$\bigcup_{i\in I}\mathcal{Z}(\Sigma_i')\cup \mathcal{Z}(\Sigma_1)=\mathcal{Z}(\Sigma).$$

We finish the proof by showing that for $\lambda\in k^q\setminus \mathcal{Z}(\Sigma)$, the polynomial
$f(\lambda,\bfx)$ has no repeated irreducible factor. Suppose it has such a factor $P\in k[\bfx]$, then choose $1\leq i\leq r$ such that $\deg_{x_i}(P)>0$. Then we must have $i\in I$. However, our choice of $\lambda$ yields that $\Res_{x_i}(f(\lambda,\bfx),f_i(\lambda,\bfx))\neq 0$ contradicting the existence of the repeated factor $P$.
\end{proof}

\section{Proof of Theorem \ref{main_thm}  }\label{Theorem1general}
\subsection{A key technical result}
First, we recall  the  following  GCD theorem.
\begin{theorem}[{\cite[Theorem 8]{GSW}}] \label{movinggcdunit}
Let $S\subset  C(\mathbf k) $ be a finite set of points. Let $F,\,G\in K[x_{1},\dots,x_{n}]$
be a coprime pair of nonconstant polynomials. For any $\epsilon>0$,
there exist an integer $m$, positive reals $c_{i}$, $1\le i\le4$,
all depending only on $\epsilon$, $n$, $\deg(F)$ and $\deg(G)$ such that for all $n$-tuple $(u_{1},\hdots,u_{n})\in({\cal O}_{S}^{*})^{n}$
with 
\begin{align}\label{multiheight10}
\max_{1\le i\le n}h(u_{i})\ge c_{1}(\tilde{h}(F)+\tilde{h}(G))+c_{2}\chi^+_S(C),
\end{align}
we have that either 
\begin{align}
h(u_{1}^{m_{1}}\cdots u_{n}^{m_{n}})\le c_{3}(\tilde{h}(F)+\tilde{h}(G))+c_{4}\chi^+_S(C)\label{multiheight11}
\end{align}
holds for an  $n$-tuple of  integers $(m_{1},\hdots,m_{n})\ne (0,\hdots,0)$   with $\sum|m_{i}|\le m$,
or
\begin{align} \label{enu:;Nsgcd<}  N_{S,{\rm gcd}}(F(u_{1},\hdots,u_{n}),G(u_{1},\hdots,u_{n}))\le\epsilon\max_{1\le i\le n}h(u_{i}).
\end{align}
\end{theorem}
\begin{remark*}
The proof in \cite[Theorem 8]{GSW} is formulated for $n\ge 2$.  When $n=1$, it is an easy consequence of the theorem of resultants.
\end{remark*}

The following key technical result will be used repeatedly in the proof of Theorem~\ref{main_thm}.
\begin{theorem}\label{main_thm_1}
Let  $G \in K[ x_1,\hdots, x_n]$ be a non-constant    polynomial  with no monomial  and repeated factors.  Then for any $\epsilon>0$,   there exist  an integer $m$ and positive reals $c_i$, $1\le i\le 4$,  all depending only on $\epsilon$, $n$, and  $\deg(G)$, 
 such that for all $n$-tuples  $\mathbf{u}:=(u_1,\dots,u_n)\in ({\cal O}_{S}^*)^n$  
 we have that either
\begin{align}\label{ht}
\max_{1\le i\le n}h(u_{i})\le  c_1  \tilde h(G) +c_2 \max\{1, \chi_S(C)\},
\end{align}
\begin{align}\label{htm}
\text{or}\ h(u_{1}^{m_{1}}\cdots u_{n}^{m_{n}})\le c_{3}  \tilde{h}(G) +c_{4} \chi_S^+(C)
\end{align}
holds for an  $n$-tuple of  integers $(m_{1},\hdots,m_{n})\ne (0,\hdots,0)$  with $\sum|m_{i}|\le m$,
 or 
\begin{align}\label{multizero}
N_{S}( G(\mathbf{u}) )-N_{S}^{(1)}( G(\mathbf{u}) )\le \epsilon \max_{1\le i\le n}\{h(u_i)\}\quad\text{ if $G(\mathbf{u})\ne 0$.}
\end{align} 
 \end{theorem}

 \begin{proof} 
  Let $\mathbf{p}\in C(\mathbf{k})$ and let $t_{\p}$ be a local parameter at $\p$. 
By  \eqref{chain rule}  and \eqref{value}, we have 
\begin{align}\label{transformpar}
\frac{d}{dt_{\p}}(G(\mathbf{u}))=G(\mathbf{u})'\cdot \frac{d t}{d t_{\p}} =D_{\mathbf{u}}(G)(\mathbf{u})\cdot \frac{d t}{d t_{\p}}.
\end{align}
 Assume that $G(\mathbf{u})\ne0$.   We claim that
 \begin{equation}\label{truncated_1}
 v_{\p}^{0}(G(\mathbf{u}))-\min\{1, v_{\p}^{0}(G(\mathbf{u}))\}\le    \min\{ v_{\p}^{0}(G(\mathbf{u})), v_{\p}^{0}(D_{\mathbf{u}}(G)(\mathbf{u})) \}+v_{\p}^0(\frac{d t}{d t_{\p}}).
\end{equation}
The above equation holds trivially if $v_{\p}^{0}(G(\mathbf{u}))\le 1$.  Therefore, we only need to consider when
$v_{\p}(G(\mathbf{u}))\ge 2$.  In this case, it is clear that $v_{\p}^{0}(G(\mathbf{u}))-1\le v_{\p}^{0}(G(\mathbf{u}))+v_{\p}^0(\frac{d t}{d t_{\p}})$ and by \eqref{transformpar} we have
$$
v_{\p}^{0}(G(\mathbf{u}))-1=v_{\p}^{0}(\frac{d}{dt_{\p}}(G(\mathbf{u}))) \le v_{\p}^{0}(D_{\mathbf{u}}(G)(\mathbf{u}))+v_{\p}^0(\frac{d t}{d t_{\p}}).
$$
Therefore, we have \eqref{truncated_1}.
Then it yields
\begin{equation}\label{trun_thm_eq_1}
    \begin{split}
        N_{S}(G(\mathbf{u}))-N_{S}^{(1)}(G(\mathbf{u}))        &\le \sum_{\p\in C({\bf k})\setminus S}\min\{ v_{\p}^{0}(G(\mathbf{u})), v_{\p}^{0}(D_{\mathbf{u}}(G)(\mathbf{u})) \}+\sum_{\p\in C({\bf k})\setminus S}v_{\p}^{0}(\frac{d t}{d t_{\p}})\\
        &\le N_{S,\gcd} (G(\mathbf{u}), D_{\mathbf{u}}(G)(\mathbf{u}))+3\gen,
    \end{split}
\end{equation}
by Proposition \ref{functiont}.   
By Lemma \ref{lem:coprime-gen},  we may factor $G=A \cdot B\in K[x_{0},\hdots,x_{n}]$ such that  $B(\mathbf{u})\in \mathcal O_S^{*}\cup\{0\}$, $A$ is constant  or  $A$ and $D_{\mathbf{u}}(G)$  
are coprime in $K[x_{1},\hdots,x_{n}]$,  and $\tilde h(A)\le 2\tilde h(G)$.  If $A$ is constant, then 
$N_{S } (G(\mathbf{u}))=N_{S } (A)\le h(A)\le 2\tilde h(G)$.  We can conclude \eqref{multizero} by assuming $h(\mathbf{u})\ge  2\epsilon^{-1}  \tilde h(G)$.  Therefore, we only need to consider when $A$ is not constant.   Then 
\begin{align}\label{trun_thm_eq_DA}
N_{S,\gcd} (G(\mathbf{u}), D_{\mathbf{u}}(G)(\mathbf{u}))= N_{S,\gcd} (A(\mathbf{u}), D_{\mathbf{u}}(G)(\mathbf{u})).
\end{align}
Since $A$ and $D_{\mathbf{u}}(G)$ are coprime, we can verify the inequality \eqref{multizero} from Theorem \ref{movinggcdunit} by potentially increasing $c_1$ and $c_2$, as supported 
  by Proposition \ref{heightDu} and the inequality $\tilde h(A)\le 2 \tilde h(G)$.
   \end{proof}

\subsection{Proof of Theorem~\ref{main_thm}}\label{generalZ123}
\begin{notation}
For a matrix $A=(a_{ij})$ with complex-valued entries, let 
$$\Vert A\Vert_{\infty}=\max_{i}\sum_{j}\vert a_{ij}\vert$$
be the maximum of the absolute row sums.
\end{notation}

The following is a simple example to illustrate the rough ideal how the exceptional sets in Theorem \ref{main_thm} are constructed: we need Theorem~\ref{main_thm_1} for the change of variables from $(x_1,x_2,x_3)$ to $(\Lambda_1,\Lambda_2,T)$ and the use of the resultant is from the proof of Lemma~\ref{special}.
 \begin{example}
We consider the example in Section~\ref{sec:intro} with $n=3$ and $G = 1 + 2x_1 - x_2 + x_3$. It is evident that the points $(u_1, u_2, u_3) = (t^\ell, at^{2\ell}, (1+a)t^{2\ell})$, $\ell\ge 1$, lies within the Zariski closed subset of $\mathbb{A}^3$ defined by the equation $x_1^2 + x_2 - x_3 = 0.$
It is clear that  $u_1^{-2}u_2=a$ and $u_2^{-1}u_3=\frac{a}{1+a}.$
  Let 
\begin{align}\label{transformE1}
\Lambda_1:=x_1^{-2}  x_2,\quad 
\Lambda_2:=x_2^{-1} x_3,\quad 
T :=x_1. 
\end{align}
Then it is easy to see that 
\begin{align}\label{transformE2}
x_1= T,\quad  x_2=\Lambda_1T^2,\quad x_3=\Lambda_1\Lambda_2T^2.
\end{align}
Using \eqref{transformE1}, we have
$$
G(T,\Lambda_1T^2,\Lambda_1\Lambda_2T^2)=1+2T-\Lambda_1T^2+\Lambda_1\Lambda_2T^2.
$$
Let $B=1+2T-\Lambda_1T^2+\Lambda_1\Lambda_2T^2\in K[\Lambda_1,\Lambda_2][T]$ be a polynomial in $T$ over $K[\Lambda_1,\Lambda_2]$. Then its discriminant (or resultant of $B$ and its derivative w.r.t. $T$) is some constant  multiple of 
$$
1+\Lambda_1-\Lambda_1\Lambda_2=1+x_1^{-2}  x_2-x_1^{-2}  x_3=x_1^{-2}(x_1^2+x_2-x_3),
$$
by \eqref{transformE2}.
\end{example}
 
\begin{proof}[Proof of Theorem~\ref{main_thm}]
When $G$ is a polynomial in $K[x]$, the result follows directly from Theorem \ref{main_thm_1}.
Now, let's consider the case where $G$ belongs to $K[x_1, \ldots, x_n]$, with $n \ge 2$. 
Let $(u_1,\ldots,u_n)\in (\mathcal{O}_S^*)^n$ such that
\begin{equation}\label{eq:c0}
\max_{1\leq i\leq n}\{h(u_i)\} > c_0 \left(\tilde h(  G) +\max\{1, \chi_S(C)\}\right)
\end{equation}
where $c_0$ is a large positive real number that will be specified later.

We will construct a proper Zariski closed subset $\tilde{Z}$ of $\mathbb{A}^n(K)$ such that (a) in part (ii) holds for $(u_1,\ldots,u_n)\notin \tilde{Z}$. Equivalently, suppose that $G(u_1,\ldots,u_n)\neq 0$ and
\begin{equation}\label{eq:(ii)(a) fails}
N_S(G(u_1,\ldots,u_n))-N_S^{(1)}(G(u_1,\ldots,u_n))>\epsilon \max_{1\leq i\leq n}\{h(u_i)\}
\end{equation}
then we show that $(u_1,\ldots,u_n)$ belongs to the desired $\tilde{Z}$. 

Our arguments are carried out inductively in several steps. In the following, the $c_{i,j}$'s and $M_i$'s denote positive real numbers depending only on $\epsilon$, $n$, $\deg(G)$, and the previously defined $c_{i',j'}$ and $M_{i'}$. It is important to note that they are independent of the (not yet specified) $c_0$.

\textbf{Step 1}: we apply Theorem~\ref{main_thm_1} to get an $n$-tuple of integers $(m_1,\ldots,m_n)\neq (0,\ldots,0)$ with $\sum\vert m_i\vert\leq M_1$ such that
\begin{equation}\label{eq:c_1,3 and c_1,4}
\lambda_1:=u_1^{m_1}\cdots u_n^{m_n}\quad\text{and}\quad h(\lambda_1)\leq c_{1,3}\tilde{h}(G)+c_{1,4}\max\{1,\chi_S(C)\};
\end{equation}
here we replace the term $\chi_S^{+}(C)$ in the statement of Theorem~\ref{main_thm_1} by the greater or equal term $\max\{1,\chi_{S}(C)\}$. 
We may assume $\gcd(m_1,\ldots,m_n)=1$. By Lemma \ref{main_lemma}, $(m_1,\hdots,m_n)$ extends to a basis 
$(m_1,\hdots,m_n)$,  $(a_{21},\hdots,a_{2n}),\hdots, (a_{n,1},\hdots,a_{n,n})$
of  $\mathbb Z^n$ such that 
\begin{align}\label{basisbound_step1}
|a_{i1}|+\cdots+|a_{in}|\le M_1+n\quad\text{ for } 2\le i\le n.
\end{align} 

Consider the change of variables
\begin{align}\label{transform_var_step1}
\Lambda_1:=x_1^{m_1}\cdots x_n^{m_n}, \quad\text{and } \quad
X_{1,i}:=x_1^{a_{i1}}\cdots x_n^{a_{in}}  \quad\text{for } 2\le i\le n
\end{align}
and put
\begin{align}\label{transform_unit_step1}
\beta_{1,i}=u_1^{a_{i1}}\cdots u_n^{a_{in}}  \quad\text{for } 2\le i\le n.
\end{align}

Let $A_1$ denote the $n\times n$ matrix whose rows are the above basis of $\mathbb{Z}^n$. Then we formally express the above identities as:
\begin{equation}\label{eq:formal_A1}
(\Lambda_1,X_{1,2},\ldots,X_{1,n})=(x_1,\ldots,x_n)^{A_1}\quad \text{and} \quad (\lambda_1,\beta_{1,2},\ldots,\beta_{1,n})=(u_1,\ldots,u_n)^{A_1}.
\end{equation}
Let $B_1=A_1^{-1}$. The entries of $B_1$ can be bounded from above in terms of $M_1$ and $n$. We have:
\begin{equation}\label{eq:formal_B1}
(x_1,\ldots,x_n)=(\Lambda_1,X_{1,2},\ldots,X_{1,n})^{B_1}\quad \text{and} \quad (u_1,\ldots,u_n)=(\lambda_1,\beta_{1,2},\ldots,\beta_{1,n})^{B_1}.
\end{equation}

Let $G_1(\Lambda_1,X_{1,2},\ldots,X_{1,n})\in K[\Lambda_1,X_{1,2},\ldots,X_{1,n}]$ with no monomial factors
and 
\begin{equation}\label{eq:G_1}
G((\Lambda_1,X_{1,2},\ldots,X_{1,n})^{B_1})=\Lambda_1^{d_1}X_{1,2}^{d_2}\cdots X_{1,n}^{d_n} G_1(\Lambda_1,X_{1,2},\ldots,X_{1,n})
\end{equation}
for some integers $d_i$, $1\leq i\leq n$. Since the transformations in \eqref{eq:formal_A1} and \eqref{eq:formal_B1} are invertible of each other and $G$ has no repeated irreducible factors, we have that $G_1$ has no repeated irreducible factors either. The coefficients of $G_1$ are the same as the coefficients of $G$ and $\deg(G_1)$ can be bounded from above explicitly in terms of $M_1$, $n$, and $\deg(G)$. Consider $G_1(\lambda_1,X_{1,2},\ldots,X_{1,n})\in K[X_{1,2},\ldots,X_{1,n}]$, by using \eqref{eq:c_1,3 and c_1,4} we have:
\begin{equation}\label{eq:tilde h after specialization 1}
\tilde{h}(G_1(\lambda_1,X_{1,2},\ldots,X_{1,n}))\leq c_{1,5}\tilde{h}(G)+c_{1,6}\max\{1,\chi_S(C)\}.
\end{equation}

For the particular change of variables in \eqref{eq:formal_A1}, \eqref{eq:formal_B1}, and \eqref{eq:G_1} (that depends on the matrix $A_1$),
we apply the lemmas in Section~\ref{mainlemmas} to obtain a proper Zariski closed subset $Z_1'$
 of
$\mathbb{A}^n(K)$ such that if $(u_1,\ldots,u_n)\notin Z_1'$ then $G_1(\lambda_1,X_{1,2},\ldots,X_{1,n})$
has neither monomial nor repeated irreducible factors. 
We now define $Z_1$ to be the union of all
such $Z_1'$ where $A_1$ ranges over the finitely many elements of $\GL_n(\mathbb{Z})$ with $\Vert A_1\Vert_{\infty}\leq M_1+n$. From the lemmas in Section~\ref{mainlemmas} and their proofs, we have that $Z_1$ satisfies properties {\rm(Z1)}--{\rm(Z3)}.

Since the $u_i$'s, $\lambda_1$, and $\beta_{1,j}$'s are in $\mathcal{O}_S^{*}$, we have:
\begin{align}\label{eq:N_S-N_S^1 step1}
\begin{split}
& N_S(G(u_1,\ldots,u_n))-N_S^{(1)}(G(u_1,\ldots,u_n))\\
=& N_S(G_1(\lambda_1,\beta_{1,2},\ldots,\beta_{1,n}))-N_S^{(1)}(G_1(\lambda_1,\beta_{1,2},\ldots,\beta_{1,n})).
\end{split}
\end{align}

From  \eqref{eq:formal_A1} and \eqref{eq:formal_B1}, we have:
\begin{equation}\label{eq:c_1,78910}
c_{1,7}\max_{1\leq i\leq n}\{h(u_i)\}-c_{1,8} \leq \max\{h(\lambda_1),h(\beta_{1,2}),\ldots,h(\beta_{1,n})\}\leq c_{1,9}\max_{1\leq i\leq n}\{h(u_i)\}+c_{1,10}.
\end{equation} 

By using \eqref{eq:c0}, \eqref{eq:c_1,3 and c_1,4}, and \eqref{eq:c_1,78910} and requiring $c_0$ be sufficiently large, we have:
\begin{equation}\label{eq:c_1,11 and c_1,12}
c_{1,11}\max_{1\leq i\leq n}\{h(u_i)\}\leq \max_{2\leq i\leq n} \{h(\beta_{1,i})\}\leq c_{1,12}\max_{1\leq i\leq n} \{h(u_i)\}.
\end{equation} 

We have the following inequalities at the end of this step:
\begin{align}\label{eq:end step1 1}
\begin{split}
	\max_{2\leq i\leq n}\{h(\beta_{1,i})\}&>c_{1,11}c_0\left(\tilde{h}(G)+\max\{1,\chi_{S}(C)\}\right)\\
										&>c_{1,13}c_0\left(\tilde{h}(G_1(\lambda_1,X_{1,2},\ldots,X_{1,n}))+\max\{1,\chi_S(C)\}\right)	
\end{split}
\end{align}
and
\begin{align}\label{eq:end step1 2}
\begin{split}
N_S(G_1(\lambda_1,\beta_{1,2},\ldots,\beta_{1,n}))-N_S^{(1)}(G_1(\lambda_1,\beta_{1,2},\ldots,\beta_{1,n}))&>\epsilon \max_{1\leq i\leq n}\{h(u_i)\}\\
&\geq \frac{\epsilon}{c_{1,12}}\max_{2\leq i\leq n}\{h(\beta_{1,i})\}.
\end{split}
\end{align}

There are $n-1$ many steps in total. Hence if $n\geq 3$, we proceed with the following $n-2$ many more steps.

\textbf{Step 2}: we include this step in order to illustrate the transition from Step $s-1$ to Step $s$
below. Since the various estimates and constructions are similar to those in Step 1, we skip some of the details. Suppose $(u_1,\ldots,u_n)\notin Z_1$ so that $G_1(\lambda_1,X_{1,2},\ldots,X_{1,n})$ has
 neither monomial nor repeated factors. 
 
 By requiring $c_0$ be sufficiently large, we can apply
 Theorem~\ref{main_thm_1} for $G_1(\lambda_1,X_{1,2},\ldots,X_{1,n})$ and $(\beta_{1,2},\ldots,\beta_{1,n})$
 and use \eqref{eq:tilde h after specialization 1}, \eqref{eq:end step1 1}, and \eqref{eq:end step1 2} to get an $(n-1)$-tuple
 $(m_2',\ldots,m_n')\neq (0,\ldots,0)$ with $\sum \vert m_i'\vert\leq M_2$ such that
  \begin{equation}\label{eq:c_2,3 and c_2,4}
\lambda_2:=\beta_{1,2}^{m_2'}\cdots \beta_{1,n}^{m_n'}\quad\text{and}\quad h(\lambda_2)\leq c_{2,3}\tilde{h}(G)+c_{2,4}\max\{1,\chi_S(C)\}.
\end{equation}
 We may assume $\gcd(m_2',\ldots,m_n')=1$. By Lemma \ref{main_lemma}, $(m_2',\hdots,m_n')$ extends to a basis of $\mathbb{Z}^{n-1}$ in which each vector has $\ell_1$-norm at most $M_2+n$.

Let $A_2'$ be the $(n-1)\times (n-1)$ matrix whose rows are the above basis of $\mathbb{Z}^{n-1}$. We make
the transformation:
$$(\Lambda_2,X_{2,3},\ldots,X_{2,n})=(X_{1,2},\ldots,X_{1,n})^{A_2'} \quad \text{and} \quad
(\lambda_2,\beta_{2,3},\ldots,\beta_{2,n})=(\beta_{1,2},\ldots,\beta_{1,n})^{A_2'}.$$
Let $A_2=(1)\oplus A_2'$ be the $n\times n$ block diagonal matrix with the $(1,1)$-entry $1$ and the matrix $A_2'$ in the remaining $(n-1)\times (n-1)$ block. We have
$$(\Lambda_1,\Lambda_2,\ldots,X_{2,n})=(\Lambda_1,X_{1,2},\ldots,X_{1,n})^{A_2} \quad \text{and} \quad
(\lambda_1,\lambda_2,\ldots,\beta_{2,n})=(\lambda_1,\beta_{1,2},\ldots,\beta_{1,n})^{A_2}.$$
Combining this with \eqref{eq:formal_A1}, we have:
\begin{align}\label{eq:formal_A2}
\begin{split}
(\Lambda_1,\Lambda_2,X_{2,3},\ldots,X_{2,n})=(x_1,\ldots,x_n)^{A_2A_1}\quad \text{and}\\
(\lambda_1,\lambda_2,\beta_{2,3},\ldots,\beta_{2,n})=(u_1,\ldots,u_n)^{A_2A_1}.
\end{split}
\end{align}
Let $B_2=(A_2A_1)^{-1}$. Let $G_2(\Lambda_1,\Lambda_2,X_{2,3},\ldots,X_{2,n})$ be the polynomial with no monomial factors such that
$$G((\Lambda_1,\Lambda_2,X_{2,3},\ldots,X_{2,n})^{B_2})=\Lambda_1^{d_1'}\Lambda_2^{d_2'}X_{2,3}^{d_3'}\cdots X_{2,n}^{d_n'}G_2(\Lambda_1,\Lambda_2,X_{2,3},\ldots,X_{2,n})$$
for some $d_1',\ldots,d_n'\in\mathbb{Z}$. We have that $\deg(G_2)$ can be bounded from above explicitly in terms of $M_2$, $M_1$, $n$, and $\deg(G)$. As before, we regard $G_2(\lambda_1,\lambda_2,X_{2,3},\ldots,X_{2,n})$ as a polynomial in $X_{2,3},\ldots,X_{2,n}$ and have the estimate:
$$\tilde{h}(G_2(\lambda_1,\lambda_2,X_{2,3},\ldots,X_{2,n}))\leq c_{2,5}\tilde{h}(G)+c_{2,6}\max\{1,\chi_S(C)\}$$
using \eqref{eq:c_1,3 and c_1,4} and \eqref{eq:c_2,3 and c_2,4}.

 For a particular $A_1$ and $A_2$, we apply the lemmas in Section~\ref{mainlemmas} to get
 a proper Zariski closed subset $Z_2'$ of $\mathbb{A}^n(K)$ such that
 if $(u_1,\ldots,u_n)\notin Z_2'$ then $G_2(\lambda_1,\lambda_2,X_{2,3},\ldots,X_{2,n})$ has neither monomial nor repeated factors. We now define $Z_2$ to be the union of all such $Z_2'$ where $A_1$ and $A_2$
 range over the finitely many unimodular matrices with
 $\Vert A_1\Vert_{\infty}\leq M_1+n$ and $\Vert A_2\Vert_{\infty}\leq M_2+n$. From the lemmas in Section~\ref{mainlemmas} and their proofs, we have that $Z_2$ satisfies properties {\rm(Z1)}--{\rm(Z3)}.
 
 By using similar estimates and requiring $c_0$ be sufficiently large as in Step 1, at the end of this step, we have
 \begin{align}\label{eq:end step2 1}
\begin{split}
	\max_{3\leq i\leq n}\{h(\beta_{2,i})\}&>c_{2,11}c_0\left(\tilde{h}(G)+\max\{1,\chi_{S}(C)\}\right)\\
										&>c_{2,13}c_0\left(\tilde{h}(G_2(\lambda_1,\lambda_2,X_{2,3},\ldots,X_{2,n}))+\max\{1,\chi_S(C)\}\right)	
\end{split}
\end{align}
and
\begin{align}\label{eq:end step2 2}
\begin{split}
N_S(G_2(\lambda_1,\lambda_2,\beta_{2,3},\ldots,\beta_{2,n}))-N_S^{(1)}(G_2(\lambda_1,\lambda_2,\beta_{2,3},\ldots,\beta_{2,n}))&>\epsilon \max_{1\leq i\leq n}\{h(u_i)\}\\
&\geq \frac{\epsilon}{c_{2,12}}\max_{3\leq i\leq n}\{h(\beta_{2,i})\}.
\end{split}
\end{align}

Let $2\leq s\leq n-1$ and suppose that we have completed Step $s-1$. This includes the construction
of $Z_{s-1}$ satisfying {\rm(Z1)}--{\rm(Z3)}. We then complete Step $s$ in the same manner Step $2$ is carried out after Step 1. The last one is Step $n-1$ resulting in the proper Zariski closed subset
$Z_{n-1}$ of $\mathbb{A}^n(K)$ satisfying {\rm(Z1)}--{\rm(Z3)}. We now define
$$\tilde{Z}=Z_1\cup \cdots\cup Z_{n-1}$$
which satisfies {\rm(Z1)}--{\rm(Z3)} since each $Z_i$ does so.
Suppose $(u_1,\ldots,u_n)\in \mathcal{O}_{S}^*\setminus \tilde{Z}$ satisfies \eqref{eq:c0}. Assume that $G(u_1,\ldots,u_n)\neq 0$ and we need to establish the inequality in part (a) of (ii). 

Assume otherwise (i.e. \eqref{eq:(ii)(a) fails} holds), then we go through all the above $(n-1)$ steps
to get the polynomial
$$P(X_{n-1,n}):=G_{n-1}(\lambda_1,\ldots,\lambda_{n-1},X_{n-1,n})\in K[X_{n-1,n}]$$ 
such that its degree can be bounded explicitly in terms of $M_{n-1},\ldots,M_1$, $n$, and $\deg(G)$. We also have
$$\tilde{h}(P)\leq c_{n-1,5}\tilde{h}(G)+c_{n-1,6}\max\{1,\chi_S(C)\}.$$
At the end of Step $n-1$, we have $\beta_{n-1,n}\in\mathcal{O}_S^*$ satisfying
 \begin{align}\label{eq:end step n-1 1}
\begin{split}
   h(\beta_{n-1,n})&>c_{n-1,11}c_0\left(\tilde{h}(G)+\max\{1,\chi_{S}(C)\}\right)\\
										&>c_{n-1,13}c_0\left(\tilde{h}(P)+\max\{1,\chi_S(C)\}\right)	
\end{split}
\end{align}
and
\begin{align}\label{eq:end step n-1 2}
N_S(P(\beta_{n-1,n}))-N_S^{(1)}(P(\beta_{n-1,n}))
> \epsilon \max_{1\leq i\leq n}\{h(u_i)\}
\geq \frac{\epsilon}{c_{n-1,12}}h(\beta_{n-1,n}).
\end{align}

Since $(u_1,\ldots,u_n)\notin Z_{n-1}$, the polynomial $P(X_{n-1,n})$
has neither monomial nor repeated irreducible factors. Let $c_{n,i}$ for $1\leq i\leq 4$ be the resulting constants when we apply Theorem~\ref{main_thm_1} for $P(X_{n-1,n})$
and $\epsilon/c_{n-1,12}$. According to Theorem~\ref{main_thm_1}, we have either
\begin{align}\label{eq:cn1 and cn2}
\begin{split}
h(\beta_{n-1,n})&\leq c_{n,1}\tilde{h}(P)+c_{n,2}\max\{1,\chi_S(C)\}\\
&\leq c_{n,1}c_{n-1,5}\tilde{h}(G)+(c_{n,1}c_{n-1,6}+c_{n,2})\max\{1,\chi_S(C)\}
\end{split}
\end{align}
or there exists a non-zero integer $r$ such that
\begin{align}\label{eq:cn3 and cn4}
	h(\beta_{n-1,n})\leq h(\beta_{n-1,n}^r) \leq c_{n,3}\tilde{h}(G)+c_{n,4}\chi_S^{+}(C).
\end{align}
We require $c_0$ be sufficiently large so that \eqref{eq:end step n-1 1} contradicts both \eqref{eq:cn1 and cn2} and \eqref{eq:cn3 and cn4}. Therefore the inequality (a) in part (ii) holds.

We now prove that the inequality (b) in part (ii) holds (after possibly enlarging $\tilde{Z}$ and increasing $c_0$ further). Write $G=\sum_{\mathbf i\in I_G}\alpha_{\mathbf i}{\mathbf x}^{\mathbf i}\in K[ x_1,\hdots, x_n]$ as in Section~\ref{sec:Preliminaries} and let $W$ be the Zariski closed subset that is the union of hypersurfaces of $\mathbb A^n$ of the form  $\sum_{\mathbf i\in J}\alpha_{\mathbf i}{\mathbf x}^{\mathbf i}=0$ where $J$ is a non-empty subset of $I_G$.
To show (ii)(b), we assume that $G(0,\hdots,0)\ne 0$ and we may further assume that $G(0,\hdots,0)=1$, which does not change $h(G)$ and may at most double $\tilde h(G)$.  
 The Zariski closed set $\tilde{Z}\cup W$ satisfies {\rm(Z1)}--{\rm(Z3)} since both $\tilde{Z}$ and $W$ do so.
Write ${\mathbf u}=(u_1,\hdots,u_n)$. Together with the condition that $\deg_{X_i}G=\deg G=d$, we may write 
$$
G( {\mathbf u})=1+\sum_{1\le i\le n}\alpha_{\mathbf i_i}u_i^d+\sum_{\mathbf i\in I_G\setminus I }\alpha_{\mathbf i}{\mathbf u}^{\mathbf i},
$$
where $I=\{(0,\hdots,0), \mathbf i_1:=(d,0,\hdots,0),\hdots, \mathbf i_n:=(0,\hdots,0,d)\}$.
 Let $S'=S\cup\{\p\notin S\,|\,v_{\p}(\alpha_{\mathbf i})\ne 0, \mathbf i\in I_G\}$.  Then  $G( {\mathbf u})\in\mathcal O_{S'}$  and
\begin{align}\label{SS}
|S'|\le |S|+2\sum_{{\mathbf i}\in I_G}  h(\alpha_{\mathbf i}).
\end{align}
For  $ {\mathbf u}\in (\mathcal O_S^*)^n\setminus (\tilde{Z}\cup W)$, we may use Theorem~\ref{BrMa}   by replacing $S$ by $S'$  to conclude that
\begin{align}\label{usingBrM}
 h(1,\alpha_{\mathbf i_1}u_1^d,\cdots,\alpha_{\mathbf i_n}u_n^d) \le N_S(G( {\mathbf u}))  +\frac{M(M-1)}{2}\max\{0,\chi_S(C)+2\sum_{{\mathbf i}\in I_G}  h(\alpha_{\mathbf i})\},
\end{align}
where $M:=\binom{n+d}{n}\ge |I_G|$.  As $h(\alpha_{\mathbf i})\le \tilde{h}(G)$, we have that 
$$
d h(1,u_1,\cdots,u_n)-\tilde{h}(G)\le N_S(G( {\mathbf u})) +\frac{M(M-1)}{2}\max\{0,\chi_S(C)+2M \tilde{h}(G) \}.
$$
Hence,
$$
N_S(G( {\mathbf u}))\ge dh(1,u_1,\cdots,u_n)- (M^3-M^2-1)  \tilde{h}(G) -\frac{M(M-1)}{2}\chi^+_S(C).
$$
We may  increase the size of $c_0$ in \eqref{eq:c0} such that  
$$
 (M^3-M^2-1) \tilde{h}(G)+\frac{M(M-1)}{2}\chi^+_S(C)\le \epsilon \max_{1\le i\le n}\{h(u_i)\}.
$$
Together with (ii)(a) and that $\max_{1\le i\le n}\{h(u_i)\}\le h(1,u_1,\cdots,u_n)$, we arrive at $N_S^{(1)}(G( {\mathbf u}))\ge (d-2\epsilon)h(1,u_1,\cdots,u_n)$. By letting $\tilde{Z}\cup W$ be the desired exceptional set $Z$, we finish the proof of the theorem without property {\rm(Z4)}.

It remains to consider the case $n=2$ and construct a slightly different exceptional set $Z$ so that properties {\rm(Z1)}--{\rm(Z4)} hold in this case. Since $n=2$, we have only Step 1 in the above arguments. To simplify the notation, we use $(\Lambda,T)$ for the indeterminates $(\Lambda_1,X_{1,2})$ in Step 1.

First, we prove that $\tilde{Z}=Z_1$ satisfies 
{\rm(Z4)}. Since $Z_1$ is the union of all the $Z_1'$ where $A_1$ ranges over all the unimodular
matrices $\begin{pmatrix} m_1&m_2 \\
a_{21} &a_{22}\end{pmatrix}$ with $\max\{\vert m_1\vert+\vert m_2\vert,\vert a_{21}\vert+\vert a_{22}\vert\}\leq M_1+2$, it suffices to fix such an $A_1$ and prove that $Z_1'$ satisfies {\rm(Z4)}. With
this $A_1$, $B_1=A_1^{-1}$, and the transformations \eqref{eq:formal_A1}--\eqref{eq:formal_B1}, we obtain the resulting $G_1(\Lambda,T)$ as in \eqref{eq:G_1}. By Lemma~\ref{special} and its proof, we obtain a proper Zariski closed subset (i.e.~a finite subset) $\mathcal{L}$ of $\mathbb{A}^1(K)$ such that $G_1(\lambda,T)$, for $\lambda\in K\setminus \mathcal{L}$, has neither monomial nor repeated factors. Then let $Z_1'$ be the union of the hypersurfaces $x_1^{m_1}x_2^{m_2}=\lambda$ for
$\lambda\in \mathcal{L}\setminus\{0\}$ (it suffices to consider $\lambda\neq 0$ since every $(u_1,u_2)\in(\mathcal{O}_S^*)^2$ automatically satisfies $u_1^{m_1}u_2^{m_2}\neq 0$). Therefore $Z_1'$ satisfies {\rm(Z4)} and so does $\tilde{Z}$.

Second, we define another proper Zariski closed subset of $\mathbb{A}^2(K)$ satisfying 
{\rm(Z1)}--{\rm(Z4)}. Fix $A_1=\begin{pmatrix} m_1&m_2 \\
a_{21} &a_{22}\end{pmatrix}$ and get $G_1(\Lambda,T)$ as before. Write 
$\displaystyle G_1=\sum_{(i,j)\in I_{G_1}}\gamma_{ij}\Lambda^iT^j$
where $\gamma_{ij}\neq 0$ for $(i,j)\in I_{G_1}$. Let $j\geq 0$ such that the set
$$I_{G_1}(j):=\{i\geq 0:\ (i,j)\in I_{G_1}\}$$
is non-empty. For each non-empty subset $I$ of $I_{G_1}(j)$, we get the non-zero polynomial:
$$P_{j,I}(\Lambda)=\sum_{i\in I} \gamma_{ij}\Lambda^i.$$
Let $\mathcal{M}$ be the finite subset of $\mathbb{A}^1(K)$ that is the union of the zero loci of all such polynomials $P_{j,I}$. Let $V'$ be the union of the hypersurfaces 
$x_1^{m_1}x_2^{m_2}=\lambda$ with $\lambda\in \mathcal{M}\setminus\{0\}$. We now let
$V$ be the union of all such $V'$ as $A_1$ ranges over all the uniomodular matrices with $\Vert A_1\Vert\leq M_1+2$. We have that $V$ satisfies {\rm(Z1)}--{\rm(Z4)} since each $V'$ does so.

Finally, we show that we can take $\tilde{Z}\cup V$ to be the exceptional Zariski set $Z$ in the current case $n=2$. Let $(u_1,u_2)\in(\mathcal{O}_S^*)^2\setminus (\tilde{Z}\cup V)$ satisfy \eqref{eq:c0} and $G(u_1,u_2)\neq 0$. We have already proved that (ii)(a) holds for $(u_1,u_2)\notin \tilde{Z}$. It remains to prove that $(u_1,u_2)\notin W$ and hence (ii)(b) holds thanks to the earlier arguments.
We go through Step 1 again to obtain $(m_1,m_2)$ with
\begin{equation}\label{eq:recall c1,3 c1,4}
\lambda_1=u_1^{m_1}u_2^{m_2}\quad \text{and} \quad h(\lambda_1)\leq c_{1,3}\tilde{h}(G)+c_{1,4}\max\{1,\chi_S(C)\},
\end{equation}
the unimodular matrices $A_1$, $B_1=A_1^{-1}$, and the polynomial $G_1(\Lambda,T)$
in \eqref{eq:G_1}. Write $G=\sum_{\bfi\in I_G}\alpha_{\bfi}\bfx^{\bfi}$ and $G_1=\sum_{(i,j)\in I_{G_1}}\gamma_{ij}\Lambda^{i}T^j$ as before.  We recall that the coefficients
$\gamma_{ij}$'s are the coefficients $\alpha_{\bfi}$'s and $\deg(G_1)$ can be bounded explicitly in terms of $M_1$ and $\deg(G)$. We also recall \eqref{eq:end step1 1}:
\begin{equation}\label{eq:recal end step1 1}
h(\beta_{1,2})>c_{1,11}c_0\left(\tilde{h}(G)+\max\{1,\chi_S(C)\}\right).
\end{equation}

If $\bfu:=(u_1,u_2)\in W$ then we have $J\subseteq I_G$ such that the subsum $\sum_{\bfi\in J}\alpha_{\bfi}\bfu^{\bfi}$ vanishes.  Then using the transformations \eqref{eq:formal_A1}--\eqref{eq:G_1}, we have $J_1\subseteq I_{G_1}$ such that
$\sum_{(i,j)\in J_1}\gamma_{ij}\lambda_1^i\beta_{1,2}^j=0$. 
Therefore $\beta_{1,2}$ is a root of the polynomial $P(X)=\sum_{(i,j)\in J_1} \gamma_{ij}\lambda_1^iX^j$. Since $(u_1,u_2)\notin V$, the above polynomial is non-zero. Moreover, its degree is bounded from above in terms of $M_1$ and $\deg(G)$ and the height of each coefficient is at most
$c_{1,13}\tilde{h}(G)+c_{1,14}\max\{1,\chi_S(C)\}$ thanks to \eqref{eq:recall c1,3 c1,4} and the earlier observation on the coefficients of $G_1$. Therefore, by \cite[p.~57]{LangDG}, we have:
\begin{equation}\label{eq:c1,15 and c1,16}
h(\beta_{1,2})\leq c_{1,15}\tilde{h}(G)+c_{1,16}\max\{1,\chi_S(C)\}.
\end{equation}
We require $c_0$ be sufficiently large so that \eqref{eq:recal end step1 1} contradicts
\eqref{eq:c1,15 and c1,16}. This means $(u_1,u_2)\notin W$. Therefore (ii)(b) holds and we finish the proof.
  \end{proof}  

 \section{Proof of Theorem \ref{toric}}\label{Theorem2general}
 
\begin{proof}[Proof of Theorem \ref{toric}]
We recall the following  setup  of finding a natural finite open covering of $X$ from the proof of \cite[Theorem 4.4]{Levin:GCD}.
Let $\Sigma$ be the fan corresponding to the smooth projective toric variety $X$.  Then there is a finite affine  covering $\{X_{\sigma} \}$ of $X$, where $\sigma\in \Sigma$ is an $n$-dimensional smooth cone with an isomorphism $i_{\sigma}: X_{\sigma}\to\mathbb A^n$.  This isomorphism restricts to an automorphism of $\mathbb G_m^n$, where we identify $\mathbb G_m^n\subset X_{\sigma}$ naturally as a subset of $X$ and 
$\mathbb G_m^n\subset\mathbb A^n$ in the standard way such that $\mathbb A^n\setminus \mathbb G_m^n$ consists of the affine coordinate hyperplanes $\{x_i=0\}$, $i=1,\hdots,n$.  Moreover, there exists a constant $c_{\sigma,A}$ and  a proper closed subset $Z_{\sigma,A}\subset\mathbb G_m^n$, depending on $\sigma$ and  $A$ such that 
\begin{align}\label{htsigma}
h(\overline{i}_{\sigma}(P))\le c_{\sigma,A} h_A(P)
\end{align}
for all $P\in \mathbb G_m^n(K)\setminus Z_{\sigma,A}\subset X(K)$,  where $\overline{i}_{\sigma}(P):=[1:i_{\sigma}(P)]\in\mathbb P^n(K)$.

The pullback $(i_{\sigma}^{-1})^*(D|_{X_{\sigma}})$ of $D$ to $\mathbb A^n$ is defined by some nonzero polynomial $f_{\sigma}\in K[x_1,\hdots,x_n]$, which does not vanish at the origin since $D$ is in general position with the boundary of $\mathbb G_m^n$ in $X$.    Furthermore,  the difference 
$\lambda_{D,\q}(P)-v_{\q} (f_{\sigma}(i_{\sigma}(P)) (=\lambda_{D,\q}(P)-v_{\q} ((i_{\sigma}^*(f_{\sigma})(P)) )$ 
is  a $M_K$-bounded function  on every $M_K$-bounded subset of $  X_{\sigma}(K)\setminus {\rm Supp} (D)$.   
On the other hand,   $P$ is identified as an $n$-tuple of $S$-units if $P\in {\mathbb G}_m^n(\mathcal O_S)$. 
Therefore, ${\mathbb G}_m^n(\mathcal O_S)\times ({M_K\setminus S}) $ is  affine $M_K$-bounded.  Consequently,  
 $$\lambda_{D,\q}(P)-v_{\q} (f_{\sigma}(i_{\sigma}(P)) )$$ 
is an $M_K$-constant for $P\in {\mathbb G}_m^n(\mathcal O_S)\setminus {\rm Supp} (D)$ and for each $\sigma\in \Sigma$.
Therefore,
\begin{align}\label{countingD}
N_{D,S} (P)-N_{D,S}^{(1)}(P) =N_S (f_{\sigma}(i_{\sigma}(P)))-N_S^{(1)} (f_{\sigma}(i_{\sigma}(P)))+O(1) 
\end{align}
for $P\in {\mathbb G}_m^n(\mathcal O_S)\setminus {\rm Supp} (D)$ and for each $\sigma\in \Sigma$.
We apply Theorem \ref{main_thm}  to find a proper Zariski closed subset  $W_{\sigma}\subset X(K)$ containing ${\rm Supp} (D)$ for each $\sigma\in \Sigma$ such that 
\begin{align}\label{gcdsigma}
N_S (f_{\sigma}(i_{\sigma}(P)))-N_S^{(1)} (f_{\sigma}(i_{\sigma}(P))) \le  \frac{\epsilon}{c_{\sigma,A}} h(\overline{i}_{\sigma}(P)) 
\end{align}
for $P\in {\mathbb G}_m^n(\mathcal O_S)\setminus W_{\sigma}$ and for each $\sigma\in \Sigma$.
  Since the number of $X_{\sigma}$ is finite, we conclude from \eqref{htsigma}, \eqref{countingD} and \eqref{gcdsigma} that there exists a proper Zariski closed subset $Z_1$ of $X$  such that 
\begin{align}\label{multiplicity}
N_{D,S} (P)-N_{D,S}^{(1)}(P)<\epsilon h_A(P)+O(1) 
\end{align}
for $P\in {\mathbb G}_m^n(\mathcal O_S)\setminus Z_1$.
  
Next  for each $\q\in S$,  there exists a finite open cover $\{U\}$ of  $X$ such that  $U\subset X_{\sigma}$ for some $\sigma\in \Sigma$ and each $ U\times\{\q\}$ is affine $M_K$-bounded. (See the proof of \cite[Chapter 10, Proposition 1.2]{LangDG}.)
Therefore, 
\begin{align}\label{lambdaq} 
\lambda_{D,\q}(P)=v_{\q}^0(f_{\sigma}(i_{\sigma}(P))+O(1)
\end{align} 
 for all $P$ on each $U(K)\setminus {\rm Supp} (D)$ with $ U\subset X_{\sigma}$ for some $\sigma\in \Sigma$.
We now apply Corollary \ref{ProximityAffine} 
respectively to each $U$ in our choice of finite cover of $X$ (by varying $\q\in S$), to find 
 a Zariski closed subset  $Z_{\sigma,\q,U}$ containing $Z_{\sigma,A}$ from \eqref{htsigma}  of $\mathbb G_m^n$ such that
 \begin{align}\label{lambdaq2} 
\lambda_{D,\q}(P) =v_{\q}^0(f_{\sigma}(i_{\sigma}(P))+O(1)
 \le \frac{\epsilon}{c_{\sigma,A}}h(\overline{i}_{\sigma}(P))+O(1) 
 <\epsilon h_A(P)+O(1)
\end{align}
 by \eqref{htsigma},  for $\q\in S$, all $P\in   U\cap \mathbb (G_m^n(\mathcal O_S)\setminus Z_{\sigma,\q,U}) \subset X_{\sigma} $.
Since $X$ is covered by finitely many  such open sets $U$, we find 
 a proper  closed subset $Z_2$ of $X$ such that 
\begin{align}\label{proxi}
m_{D,S}(P )< |S|\epsilon  h_A(P)+O(1)
\end{align}
for all $P\in  \mathbb  G_m^n(\mathcal O_S)\setminus Z_2 $.
Since $h_D(P)= m_{D,S}(P )+N_{D,S}(P )+O(1)$, we can derive from \eqref{proxi} that
$$
N_{D,S} (P)\ge h_D(P)-|S|\epsilon h_A(P)-O(1) 
$$
for all $P\in  \mathbb  G_m^n(\mathcal O_S)\setminus Z_2$. Then \eqref{multiplicity} implies 
$$
N_{D,S}^{(1)}(P)\ge h_D(P)- (1+|S|)\epsilon h_A(P)-O(1), 
$$
for all $P\in  \mathbb  G_m^n(\mathcal O_S)\setminus \{Z_1\cup Z_2 \}$.
\end{proof}

\section{Lang-Vojta Conjecture over function fields}\label{LangVojtaF}
\subsection{Proof of Theorem \ref{thmVconj}}\label{ProofTheorem2}
While it is possible to treat Theorem \ref{thmVconj} as an application of Theorem \ref{vojtaconj}, we choose to provide a direct proof. This approach allows us to explicitly track the exceptional set and all the constants involved. Furthermore, this proof can be readily adapted to establish the corresponding Green-Griffith conjecture with an explicit exceptional set.  

\begin{proof}[Proof of Theorem \ref{thmVconj}]
  Let $\mathbf x\in \Sigma$. As noted in the paragraph following Definition \ref{int_set}, by enlarging $S$ (according to $\Sigma$), 
  we have  $N_{D_i,S}(\mathbf x)=0$ for $1\le i\le n+1$.  
  We may  assume that the coefficients of $F_i$, $1\le i\le n+1$, are in $\mathcal O_S^*$ by increasing $S$.
We note that the size of $S$  increases at most by a constant multiple of  $\sum_{i=1}^{n+1} \tilde h(F_i)$.
Let $b_{i}:= \operatorname{lcm}(\deg F_{1},\hdots,\deg F_{n+1})/\deg F_{i}$, $1\le i\le n+1$.
Then $N_{D_i,S}(\mathbf{x})=0$  implies that $ F_{i}^{b_{i}}(\mathbf{x})/F_{n+1}^{b_{n+1}}(\mathbf{x}) \in\mathcal O_S^*$ for $1\le i\le n+1$.  
Let $X=\mathbb P^n\setminus D$. 
Then we have a finite morphism
$\pi:X\to \mathbb A^n$ over $K$ given by 
$$
P\mapsto (F_{1}^{b_{1}}(P)/F_{n+1}^{b_{n+1}}(P),\hdots,F_{n}^{b_{n}}(P)/F_{n+1}^{b_{n+1}}(P)),
$$
which maps   $\mathbf{x}$ into $(\mathcal O_S^*)^n$.
Since $D_1(\p),\dots,D_{n+1}(\p)$ intersect transversally, they are in general position and hence $D_1,\hdots,D_{n+1}$ are in general position by Proposition \ref{generalposition}.
Thus the map $\pi$  extends to a finite morphism  $\tilde \pi:\PP^{n}(\overline{K})\to\PP^{n}(\overline{K})$ over $K$ given by 
$$
P\mapsto[F_{1}^{b_{1}}(P):\cdots:F_{n+1}^{b_{n+1}}(P)].
$$
Let $\mathcal Z$ be the ramification divisor of $\pi$, and $\tilde {\mathcal Z}$ be its Zariski-closure in $\mathbb P^n$.  In what follows, we will show that $\tilde {\mathcal Z}$  is nontrivial (hence ample) and  is in general position with $D$.
Moreover, we will show that $N_{\tilde Z,S}(\mathbf x)$ is ``small"    if  $\mathbf x$ is outside  $\pi^{-1}(\tilde {\mathcal Z})$.

Let $J\in K[x_0,\dots,x_n]$ be the determinant of the Jacobian matrix $\big(\frac{\partial F_i^{b_{i}}}{\partial x_j}\big)_{1\le i\le n+1, 0\le j\le n}$ of $\tilde\pi$.  Then $\tilde {\mathcal Z}$ is associated to the largest factor of $J\in K[x_0,\dots,x_n]$ that is not divisible by any of the $F_i$.  Therefore,   $\tilde {\mathcal Z}$ is contained the zero locus of  of  $ G:=\det M$, where 
\begin{equation}\label{M0}
 G:= \det   \left( \frac{\partial F_i}{\partial x_j} \right)_{1\le i\le n+1, 0\le j\le n}.
 \end{equation}

Let $\p\in C$  be as in the assumption. The specialization $D_i(\p)$ at $\p$ is the zero locus of $ F_{i,\p}(\p)$ in $\mathbb P^n({\bf k})$, where $F_{i,\p}=t_{\p}^{-v_{\p}(F_i)}  \cdot F_i$.
Let
\begin{equation}\label{Mp}
 G_{\p}:=\det\left( \frac{\partial F_{i,\p}(\p)}{\partial x_j} \right)_{1\le i\le n+1, 0\le j\le n}.
\end{equation}

Note that $\deg G_{\p}=\deg G\ge 1$  since $\sum_{i=1}^{n+1}\deg F_{i,\p}(\p)=\sum_{i=1}^{n+1}\deg F_i\ge n+2$.  By Proposition \ref{generalposition}, to show that $[G=0]$, $D_1,\hdots,D_{n+1}$ are   in general position (in $\mathbb P^n(\overline{K})$), it suffices to show that $[G_{\p}=0]$, $D_1(\p),\hdots,D_{n+1}(\p)$ are    in general position  (in $\mathbb P^n({\bf k})$).  The latter can be verified by showing that 
 $G_{\p}$ does not vanish at any intersection point of any $n$ divisors among $D_1(\p),\hdots,D_{n+1}(\p)$.  By rearranging the index, it suffices to  show that $G_{\p}(Q)\ne 0$ for  $Q\in \cap_{i=1}^n D_i(\p)$.  Since $D_1(\p),\hdots,D_{n+1}(\p)$ intersect transversally, $F_{n+1,\p}(\p)(Q)\ne 0$. 
Using the Euler formula 
$$
\sum_{j=0}^{n}\frac{\partial F_{i,\p}(\p)}{\partial x_j}x_j=\deg F_i \cdot F_{i,\p}(\p), 
$$
we obtain 
\begin{equation*}
\begin{split}
        x_0G_{\p} =\det \begin{pmatrix}
        d_1 F_{1,\p}(\p)& \frac{\partial F_{1,\p}(\p)}{\partial x_1} & \cdots & \frac{\partial F_{1,\p}(\p)}{\partial x_n} \\
        \vdots & \vdots & \ddots & \vdots\\
        d_{n+1}F_{n+1,\p}(\p) & \frac{\partial F_{n+1,\p}(\p) }{\partial x_1} & \cdots & \frac{\partial F_{n+1,\p}(\p) }{\partial x_n} 
    \end{pmatrix}
\end{split}
\end{equation*}
and hence
$$
 x_0(Q)G_{\p} (Q)=(-1)^{n+1} d_{n+1}F_{n+1,\p}(\p) (Q)\det \left( \frac{\partial F_{i,\p}(\p) }{\partial x_j}(Q) \right)_{1\le i, j\le n},
$$
where $d_i:=\deg F_i$.
Since $D_1(\p),\hdots,D_{n+1}(\p)$  intersect transversally, $\det \left( \frac{\partial F_{i,\p}(\p) }{\partial x_j}(Q) \right)_{1\le i, j\le n}\ne 0$.  Then $G_{\p} (Q)\ne 0$ as  $F_{n+1,\p}(\p) (Q)\ne 0$.  
In conclusion, $\tilde {\mathcal Z}=[G=0]$, $D_1,\hdots,D_{n+1}$ are in general position.  
 
 Let $\mathcal Z_0$ be a component of ${\mathcal Z}$ such that its Zariski-closure $\tilde {\mathcal Z_0}$ in $\mathbb P^n$ is defined by an irreducible factor $G_0$ of  $G$ in $K[x_0,\hdots,x_n]$. 
Let $A\in K[y_1,\hdots,y_n]$ be the defining equation of $\pi(\mathcal Z_0)$ in $\mathbb A^n$.  Then the multiplicity of $\pi^*A$ along $\mathcal Z_0$ is at least 2.  Therefore, we have the following expression
$\pi^* A=G_0^2 H\cdot F_{n+1}^{\ell}$, where $\ell$ is a negative integer, $H$ is a homogeneous polynomial  in $K[x_0,\hdots,x_n]$, and $2\deg G_0+\deg H+\ell\deg F_{n+1}=0$.  We may assume that the coefficients of $H$ is in $\mathcal O_S$ by multiplying a suitable scalar to $A$.  We can also assume that the coefficients of $G_0$ is in $\mathcal O_S^*$ by enlarging $S$ and similarly the increasing size is bounded by a constant multiple of $\sum_{i=1}^{n+1} \tilde h(F_i)$.
  Recalling that $\mathbf{x}\in \Sigma$, we let 
  \begin{align}\label{ui}
   u_i:=F_{i}^{b_{i}}(\mathbf{x})/F_{n+1}^{b_{n+1}}(\mathbf{x}) \in\mathcal O_S^*,
   \end{align} 
   for $1\le i\le n $, and let $\mathbf{u}=(u_1,\hdots,u_n)$.  
Then 
$$
A(\mathbf{u})= G_0^2(\mathbf{x})H(\mathbf{x})F_{n+1}^{\ell}(\mathbf{x}).
$$
Let $\q\in C({\bf k})\setminus S$.   To estimate $\lambda_{\tilde{\mathcal Z_0},\q} (\mathbf{x})$,  we may assume that $e_{\q}(\mathbf x):=\min_{0\le i\le n}\{v_{\q}(x_i)\}=0$ by multiplying each $x_i$ by $t_{\q}^{-e_{\q}(\mathbf x)}$. Since $\lambda_{D_{n+1},\q}(\mathbf x)=0$ and the coefficients of $F_{n+1}$ are in $\mathcal O_S^*$, we have $v_{\q}(F_{n+1}(\mathbf x))=0$. 
Together with the fact that the coefficients of $H$ are in $\mathcal O_S$, we see that $v_{\q}(A(\mathbf{u}))\ge 2v_{\q}( G_0(\mathbf{x}))\ge0$; this shows that 
 $\lambda_{\tilde{\mathcal Z_0},\q} (\mathbf{x})=v_{\q}(G_0(\mathbf{x}))\le v_{\q}( A(\mathbf{u}))-\min\{1,v_{\q}( A(\mathbf{u}))\}$.  Therefore,
\begin{equation}\label{ramificationcounting}
    N_{ \tilde{\mathcal Z_0}, S} (\mathbf{x})\le N_S(  A(\mathbf{u}))-N_S^{(1)}(  A(\mathbf{u})).
\end{equation}
We now apply  Theorem \ref{main_thm}  to  $A$ evaluating at $\mathbf{u}$.  Recall that $[1:u_1:\cdots:u_n]= [F_1^{b_1}(\mathbf{x}):\cdots:F_{n+1}^{b_{n+1}}(\mathbf{x})]$ in $\mathbb P^n$. Since these $F_i$ are in general position, we have
\begin{align}\label{htrelation}
bh(\mathbf{x})-c'\le h(1,u_1,\hdots,u_n)\le bh(\mathbf{x})+c',
\end{align}
where $b=\operatorname{lcm}(\deg F_{1},\hdots,\deg F_{n+1})$ and $c'$ is a positive real depending only on these $F_i$.    Then for
any $\epsilon>0$,  there exists an integer $m$, positive reals $c_1'$ and $c_2'$, and
 a proper Zariski closed subset   $Z'$ of $\mathbb P^n$ over $K$, all depending only on $\epsilon$ and those $F_i$, such that   if     $\mathbf{x}\notin Z'\cup \pi^{-1}(\pi(\mathcal Z_0))$, then together with \eqref{htrelation} we have that either 
\begin{align}\label{heightupper}
h(\mathbf{x})\le c_{2}'\max\{1,\chi_S(C)\},
\end{align}
or 
\begin{equation}\label{NS1}
  N_{ \tilde{\mathcal Z_0}, S} (\mathbf{x})\le N_{S}(A(\mathbf{u}))-N_{S}^{(1)}(A(\mathbf{u}))\le \epsilon h(\mathbf{x}), 
\end{equation}
and
\begin{equation}\label{NS2}
N_{S}^{(1)}(A(\mathbf{u}))\ge \deg A(1-\epsilon)  h (1,u_1,\hdots,u_n). 
\end{equation}
Denote by $\tilde A$ the homogenization of $A$.  By increasing $c_2'$ in \eqref{heightupper}, we may assume that the coefficients of $A$ are $S$-units.
Then by that $\tilde{\mathcal Z_0}\le \tilde \pi^*([\tilde A=0])$,   functorial  property and \eqref{NS2}, we have 
\begin{align}\label{MS3}
m_{ \tilde{\mathcal Z_0}, S} (\mathbf{x})&\le m_{ [\tilde A=0], S} (\pi(\mathbf{x}))
=\deg A\cdot h (1,u_1,\hdots,u_n)-N_{S}(A(\mathbf{u}))+h(A) \cr
& \le  \deg A\cdot\epsilon h (1,u_1,\hdots,u_n)+h(A).
\end{align}
Since $\mathcal Z_0$ is defined by $G_0$, it then follows from \eqref{NS1}, \eqref{MS3} and \eqref{htrelation} that
\begin{equation}\label{htbdd2}
\deg G_0 \cdot h (\mathbf{x})+ h(G_0)=h_{ \tilde{\mathcal Z_0}} (\mathbf{x})\le (b\deg A+1)\epsilon h(\mathbf{x})+c'\deg A\epsilon+h(A).
\end{equation}
if $\mathbf{x}\notin Z'\cup \pi^{-1}(\pi(\mathcal Z_0))$.
Since the constants in the left hand side of \eqref{htbdd2} is independent of $\epsilon$, we may reach our assertion by making the positive $\epsilon$ sufficiently small.
 \end{proof}

\subsection{Proof of Theorem \ref{vojtaconj}}\label{ProofTheorem3}

The following lemma is the general case of \cite[Lemma 1]{CZ2013}, which treats the case  $n=2$ and $X$ over $\bf k$, and the proof extends to more general situations. 
\begin{lemma}\label{big}
Let  $ X$ be a smooth variety over $K$ of dimension $n\ge 2$, $D$ be a  normal crossings divisor on $ X$, and $V= X\setminus D$.   Suppose that $V$ admits a finite morphism  $\pi: V\to \mathbb G_m^n$ over $K$, which extends to a morphism  $\tilde\pi:  X\to \mathbb P^n$.  Let ${ Z}\subset V$ be the ramification divisor of $\pi$ and $\tilde { Z}$ be the Zariski-closure of ${ Z}$ in $ X$.  Let $H_1,\hdots,H_{n+1}$ be $n+1$ hyperplanes in general position in $\mathbb P^n$, and  assume that $D$ equals the support of $ \tilde\pi^*H_1+\hdots+\tilde\pi^*H_{n+1}$ (i.e. the sum of components of $ \tilde\pi^*H_1+\hdots+\tilde\pi^*H_{n+1}$ counted with multiplicity 1).   Then $\tilde {  Z}$ is linearly equivalent to $D+ {\mathbf K}_{ X}$.  In particular, if $V$ is of log-general type, then $\tilde {  Z}$ is big.
\end{lemma}

\begin{proof}
We will follow the  arguments  from \cite{CZ2013}.
By \cite[(1.11)]{debarre2001-higher}, the canonical class on $ X$ could be written as 
\begin{equation*}
    {\mathbf K}_{ X}\sim \tilde \pi^*({\mathbf K}_{\PP^n})+	  {\rm Ram},
\end{equation*}
where $ {\rm Ram}$ is the ramification divisor of $\tilde \pi: X\to \PP^n$. Note that  $ {\rm Ram}$ decomposes as $ {\rm Ram}=\tilde{Z}+R_D$, where $R_D$ is the contribution coming from the support contained in $D$. On the other hand, 
we also have ${\mathbf K}_{\PP^n}\sim -(H_1+\cdots+H_{n+1})$ and  $\tilde{\pi}^*(H_1+\cdots+H_{n+1})=D+R_D$. Then we obtain $\tilde{Z}\sim {\mathbf K}_{ X}+D $. 
\end{proof}
 
\begin{proof}[Proof of Theorem \ref{vojtaconj}]
 Let $D$ be a  normal crossings divisor on $\ X$. Let $\tilde \pi:  X\to \mathbb P^n$  be a morphism over $K$ such that its restriction to $V:= X\setminus D$, $\pi: V\to \mathbb G_m^n$,  is  a finite morphism.
Let $H_1,\hdots,H_{n+1}$ be the $n+1$ coordinate hyperplane divisors of   $ \mathbb P^n$.  
 Since a finite morphism is a closed map, we see that $\pi$ is surjective as  $\dim V=\dim\pi(V)=\dim \mathbb G_m^n.$ 
Then  $D$ and $\tilde \pi^* (H_1+\cdots+H_{n+1})$ have the same support.
 
Recall that  $\Sigma\subseteq  X(K)$ is a  $(D,S)$-integral  subset. By Proposition  \ref{under_mor}, we see that $\tilde\pi(\Sigma)\subset\PP^n(K)$ is a  $(D',S)$-integral subset, where  $D'=H_1+\cdots+H_{n+1}$.
Let  ${\mathbf x}\in \Sigma\subset V$ and $\pi(\mathbf x)=\mathbf y:=[y_0:\cdots:y_n]$.  Then, by enlarging $S$, we have 
$$
\lambda_{D',\p}(\mathbf y)=\left(v_{\p}(y_0)-v_{\p}(\mathbf y)\right)+\cdots+\left(v_{\p}(y_n)-v_{\p}(\mathbf y)\right)=0 
$$
for $\p\notin S$ as noted in the paragraph following Definition \ref{int_set}.  As $v_{\p}(y_i)\ge v_{\p}(\mathbf y)$, 
we have $v_{\p}(y_0)=\cdots=v_{\p}(y_n)=v_{\p}(\mathbf y)$ for $\p\notin S$, and hence $v_{\p}(y_i/y_0)=0$ for $0\le i\le n$ and $\p\notin S$.  Then we   have that 
\begin{align}\label{setup}
\pi({\mathbf x})=[1:u_1:\cdots:u_n]\in\mathbb P^n(K), \quad\text{where }  u_i\in\mathcal O_S^*  \text{ for }  1\le i\le n. 
\end{align}
 
Let $  Z\subset V$ be the ramification divisor of $\pi$ and $R$ be the closure of $\pi(Z)$ in $\mathbb P^n$,  which is then defined by a   homogeneous polynomial $F \in K[x_0,x_1,\hdots,x_n]$.  
Let $\tilde { Z}$ be the Zariski-closure of $Z$ in $ X$. Then $\tilde { Z}$ is big by Lemma \ref{big}. Let 
  $\tilde {Z_0}$ be an   irreducible component of $\tilde {Z}$.  We note that since  $\tilde {Z}$ is defined over $K$, we may find a finite Galois  extension $K'$ of $K$ such that all the irreducible components of $\tilde {Z}$ are defined over $K'$.  
Then the sum over all Galois conjugate of $\tilde {Z_0}$ is a divisor over $K$.  Since  valuations are compatible with field extensions,  at each point in $  X(K)\setminus {\rm Supp} ( \tilde {Z_0})$,  we may first computing Weil functions  of $\tilde {Z_0}$ with valuation $w$ in $K'$ and then adding all the Galois conjugates of  $\tilde {Z_0}$ and $w$ together.   
Therefore, we may assume that  each irreducible component of   $\tilde{Z}$   is defined over $K$.  We note that $[F=0]$  and the coordinate hyperplanes $[x_i=0]$, $0\le i\le n$ are in general position in $\mathbb P^n(\overline K)$ by observing that $F$ does not vanish identically when evaluating at  the set of points  $\{(1,0,\hdots,0),\hdots,(0,\hdots,0,1)\} $ because of  the assumption that    the specialization $R(\p)\in \mathbb P^n({\bf k}) $  for some $\p\in C$ does not intersect the set of points  $\{(1,0,\hdots,0),\hdots,(0,\hdots,0,1)\} $.

 Suppose that ${\mathbf x}\notin{\rm Supp} (\tilde {Z_0})$.  Let ${\mathbf y}\in   Z\cap  \tilde Z_0$. 
Let  $F_0\in K[x_0,x_1,\hdots,x_n]$ be the irreducible factor of $F$ in $K[x_0,x_1,\hdots,x_n]$ such that $F_0( \pi({\mathbf y}) )=0$ and $R_0$ be the divisor on $\mathbb P^n$ defined by $F_0=0$.   Then $\tilde \pi^*R_0$ has multiplicity at least 2 along $\tilde Z_0 $.  Let $(U,f)$ be a local representation of $\tilde Z_0 $ such that  $(U,f^2g)$ is a local representation of $\tilde \pi^*R_0$, where $g$ is a regular function on $U$ over $K$.
Then both the differences
$\lambda_{\tilde Z_0,\q}({\mathbf x})-v_{\q} (f({\mathbf x}))$ and $\lambda_{\tilde \pi^*R_0,\q}({\mathbf x})-2v_{\q} (f({\mathbf x}))-v_{\p} (g({\mathbf x}))$
are $M_K$-bounded functions  on every $M_K$-bounded subset of $ X(K)\setminus {\rm Supp} (\tilde \pi^*R_0)$.
 Since  $ X(K)\times M_K$   is a finite union of affine $M_K$-bounded sets  (see 
\cite[Chapter 10, Proposition 1.2]{LangDG}), 
we have 
\begin{align}\label{twice}
0\le 2\lambda_{\tilde Z_0,\q}({\mathbf x})\le  \lambda_{\tilde \pi^*R_0,\q}({\mathbf x}) 
\end{align}
up to  $M_K$ constants.  As each $M_K$ constant has only finitely many nonzero components, by enlarging  $S$, we may  assume that (\ref{twice}) holds exactly for each $\q\in  C({\bf k})\setminus S$.
 Recall   that $\pi({\mathbf x})=[1:u_1:\cdots:u_n]\in\mathbb P^n(K)$, where $u_i\in\mathcal O_S^*$ for $1\le i\le n$. By further enlarging  $S$, we may  assume that each coefficient of $F_0$ lies in $\mathcal O_S^*$. Thus  if ${\mathbf x}\notin {\rm Supp} (\tilde \pi^*R_0)$, then for each $\q\in  C({\bf k})\setminus S$ we have    
\begin{align}\label{withF}
 \lambda_{\tilde \pi^*R_0,\q }( \mathbf x) =v_{\q}(F_0(1,u_1,\hdots,u_n)).  
\end{align}
 
 Then  \eqref{twice} implies
\begin{align}\label{coutingzero}
N_{\tilde {Z_0},S}({\mathbf x})\le N_S( F_0(1,u_1,\hdots,u_n))-N^{(1)}_{S}( F_0(1,u_1,\hdots,u_n)).
\end{align}
 Fix a positive $\epsilon$ (to be determine later). 
We now apply Theorem \ref{main_thm} to find positive constants $c_1$ and $c_2$, and   a proper Zariski closed subset $Y_1$  of $\mathbb P^n$  such that  we have  either  
\begin{align}\label{hest}
h(1,u_1,\hdots,u_n)\le   c_1  +c_2\max\{1,\chi_S(C)\},
\end{align} or 
\begin{align}\label{countingepsilon}
N_S( F_0(1,u_1,\hdots,u_n))-N^{(1)}_{S}( F_0(1,u_1,\hdots,u_n))\le \frac{\epsilon}{r}  h(1,u_1,\hdots,u_n),
\end{align}
where $r$ is the number of irreducible components of $\tilde { Z}$, and
\begin{align}\label{countingepsilon2}
N^{(1)}_{S}( F_0(1,u_1,\hdots,u_n))\ge (1-\epsilon)\deg F_0\cdot h(1,u_1,\hdots,u_n), 
\end{align}
if $\pi({\mathbf x})\notin Y_1$. 
Then by  \eqref{coutingzero} and \eqref{countingepsilon}, we have
\begin{align*}
N_{\tilde { Z_0},S}({\mathbf x})\le \frac{\epsilon}{ r}h(1,u_1,\hdots,u_n), 
\end{align*}
for ${\mathbf x} \notin \tilde { Z_0}\cup\pi^{-1}(Y_1)$.   By repeating the  above argument  for each component of $\tilde { Z}$, and by enlarging  $c_1$ and $c_2$, we find  a proper Zariski closed subset $Y$  of $\mathbb P^n$    such that   if   ${\mathbf x}\notin\tilde Z\cup\pi^{-1}(Y)$,  then  either \eqref{hest} holds or 
\begin{align}\label{countingRam}
N_{\tilde { Z},S}({\mathbf x})\le  \epsilon   h(1,u_1,\hdots,u_n), 
\end{align}
and
\begin{align}\label{countingepsilon3}
N_{S}( F(1,u_1,\hdots,u_n))\ge (1-\epsilon)\deg F\cdot h(1,u_1,\hdots,u_n). 
\end{align}
 By the functorial property,  $ \tilde { Z}\le \tilde \pi^* ([F=0])$ (as divisors)  implies that   for some reals $c_3$ and $c_4$, not depending on ${\bf x}$, we have 
\begin{align}\label{Fproxi2}
  m_{\tilde {Z},S}({\mathbf x})&\le m_{\tilde \pi^* ([F=0]),S}( \mathbf x ) +c_3 \cr
  &=m_{ [F=0],S}(1,u_1,\hdots,u_n) +c_4\cr
  &=\deg F\cdot h(1,u_1,\hdots,u_n)-N_{ [F=0],S}(1,u_1,\hdots,u_n)+\tilde h(F)+c_4\cr
  &\le \deg F\cdot\epsilon h(1,u_1,\hdots,u_n)+\tilde h(F)+c_4.  \quad\text{ (by \eqref{countingepsilon3})}
\end{align}
 Hence  if ${\mathbf x}\notin  {\rm Supp} ( \tilde\pi^*(R))\cup\pi^{-1}(Y)$,  then     either \eqref{hest} holds or
\begin{align}\label{bdd2}
h_{ \tilde {Z}}( \mathbf x ) \le (\deg F+1)\epsilon\cdot h(1,u_1,\hdots,u_n)+  \tilde h(F)+c_4 
 \end{align}
by \eqref{countingRam} and \eqref{Fproxi2}.
Since $ \tilde { Z}$ is big, there exist   positive  constants $b_1$ and $c_5$ and a   proper Zariski-closed set   $W_1$ of $ X$   over $K$, depending only on $\pi$, $F$ and $\tilde  Z$, such that 
\begin{equation}\label{bdd_ample}
    h( 1,u_1,\hdots,u_n)=\frac 1{\deg F}  \cdot h_{\tilde \pi^* ([F=0])}( \mathbf x )\le b_1 h_{\tilde{Z}}( \mathbf x)+c_5.
\end{equation} 
for all $ \mathbf x \in  X(K)$ outside  $W_1$. (See \cite[Proposition 10.11]{Vojta}.)
Then for $\epsilon<(3\deg F\cdot b_1)^{-1}$,  we deduce from \eqref{bdd2} and \eqref{bdd_ample} that
$h_{ \tilde {Z}}( \mathbf x ) \le c_6;$
and from \eqref{hest} and   \eqref{bdd2} to derive
$h_{ \tilde {Z}}( \mathbf x ) \le  c_7  +c_8\max\{1,\chi_S(C)\}$,
where $c_i$'s are positive constants independent of $\mathbf x$.
In conclusion, we have
\begin{align}\label{htZ}
h_{ \tilde {Z}}( \mathbf x ) \le  \max\{c_6,c_7\}  +c_8\max\{1,\chi_S(C)\}
\end{align}
for  ${\mathbf x}\notin W_1 \cup  \tilde Z\cup\pi^{-1}(Y)$.  Finally, let $A$ be a big  divisor on $ X$.  We can conclude our assertion  of \eqref{htboundedA} from \eqref{htZ} with the following height inequality.
By \cite[Proposition 10.11]{Vojta}, 
there exist  positive  constants $b_2 $ and $c_8$, and a   proper Zariski-closed set   $W_2$ of $ X$   over $K$, depending only on $A$ and $\tilde  Z$, such that 
\begin{equation}\label{bdd_big}
  h_A( \mathbf x )\le b_2  h_{\tilde{Z}}( \mathbf x ) +c_8
\end{equation}
for all $ \mathbf x \in  X(K)$ outside  $W_2$.  (See \cite[Proposition 10.11]{Vojta}.)
\end{proof}

\section{Remarks on the Complex Cases}\label{ComplexRemark}
 
We will now explain (without proof) the application of the reduction method in proving Theorem \ref{main_thm} in the complex case. This will enable us to derive explicit exceptional sets for both Vojta's general abc conjecture and the Green-Griffith-Lang conjecture corresponding to the split cases of Theorem  \ref{toric}-\ref{vojtaconj}.

First of all, we will need the following theorem analogous to Theorem \ref{main_thm_1}, which can be deduced by combining  the proofs of  \cite[Theorem 5.2]{GW22} and \cite[Theorem 1.1]{GSW22}. 

\begin{theorem}\label{main_thm_1G}
	Let $G$ be a non-constant  polynomial   in $\mathbb C[x_1,\hdots,x_n]$ with no monomial factors and no repeated factors.  
	Then for any   $\epsilon >0$, there exists a positive integer    $m$  such that
	for any $n$-tuple $(u_1,\dots,u_n)$  of units, i.e. entire functions with no zero, we have  either the following inequality
\begin{align}\label{exception0}
		T_{u_1^{m_1}\cdots u_n^{m_n}}\le_{\rm exc}{\rm o}(\max_{1\le i\le n}\{T_{u_i}(r)\}), 
\end{align}
	for a non-trivial $n$-tuple $(m_1,\hdots,m_n)$ of integers with  $\sum_{i=1}^n |m_i|\le  m$,  or
\begin{align}\label{multizerobdd}
	N_{G(u_1,\dots,u_n)}(0,r)-N^{(1)}_{G(u_1,\dots,u_n)}(0,r)\le_{\rm exc} \epsilon \max_{1\le i\le n}\{T_{u_i}(r)\}. 
\end{align}
\end{theorem}

We can apply the induction process employed in the proof of Theorem \ref{main_thm} (ii)(a) to deduce that the exceptional set for the aforementioned theorem is indeed a proper Zariski closed subset. Moreover, we can also adapt  the arguments presented for Theorem \ref{main_thm} (ii) (b) to the complex case. Consequently, we arrive at a   version of Conjecture \ref{ConjABC} \ref{truncate2} as follows. 

\begin{theorem}\label{main_thm_1GExc}
	Let $G$ be a non-constant  homogeneous polynomial   in $\mathbb C[x_0,\hdots,x_n]$ with no monomial factors and no repeated factors. 	Then for any   $\epsilon >0$, there exists a proper Zariski closed subset $Z$ of $\mathbb P^n$ such that for any $n+1$-tuple  $\mathbf u:=(u_0,\dots,u_n)$  of units with $\mathbf u(\mathbb C)$ not contained in $Z$, we have 
\begin{align}\label{multizerobdd1}
	N_{G(\mathbf u)}(0,r)-N^{(1)}_{G(\mathbf u)}(0,r)\le_{\rm exc} \epsilon  T_{\mathbf u}(r),   
\end{align}
and
\begin{align}\label{truncate1_2}
N^{(1)}_{G(\mathbf u)}(0,r)\ge_{\rm exc}  (\deg  G- \epsilon)\cdot T_{\mathbf u}(r),
\end{align}
 if we assume furthermore that  the hypersurface  defined by   $G$ in $\mathbb P^n$ and the coordinate hyperplanes are in general position.
\end{theorem}

The theorem stated above can be extended to encompass more general cases of toric varieties, as  Theorem \ref{toric}, though we shall omit the details here.

By using Theorem \ref{main_thm_1GExc} in the proof of \cite[Theorem 1.2]{GSW22}, we may extend the result of Noguchi, Winkelman and Yamanoi  in \cite[Theorem 5.4]{noguchi2007degeneracy} to the following  case of the (Strong) Green-Griffiths-Lang conjecture.(See \cite[Conjecture 6.6.30]{noguchibook}.)
\begin{theorem}\label{GG_conj}
Let $\mathbf{f}=(f_0,\hdots,f_n):\CC\to \PP^n$ be a holomorphic map, where $f_0,\hdots,f_n$ are entire functions without common zeros.    Let 
$F_i$, $1\le i\le n+1$ be homogeneous irreducible polynomials of positive degree  in $\mathbb C[x_0,\hdots,x_n]$ such that $\sum_{i=1}^{n+1}\deg F_i\ge n+2$. 
Assume that  the zero locus of  $F_i$,  $1\le i\le n+1$,  intersect transversally.   
Then there exists a proper Zariski closed subset  $Z$ of $\mathbb P^n$ such that the image of every entire curve
$f:\mathbb C\to \mathbb P^n(\mathbb C)\setminus \cup_{i=1}^{n+1} [F_i=0]$
is contained in $Z$.
\end{theorem}

Our proof of Theorem \ref{vojtaconj} can be adapted to establish the Green-Griffiths-Lang conjecture for varieties of log general type that are ramified covers of $\mathbb{G}_m^n$ as follows.  
 
\begin{theorem}\label{ramified}
Let  $ X$ be a complex smooth  projective  variety  of dimension $n\ge 2$, $D$ be a  normal crossings divisor on $ X$, and $V= X\setminus D$. 
Suppose there is a morphism $\tilde \pi:    X\to\mathbb P^n$ such that $\pi=\tilde \pi|_V : V\to \mathbb G_m^n$ is a finite morphism. 
Let $ Z\subset V$ be the ramification divisor of $\pi$, and $R$ be the closure of $\pi(Z)$ in $\mathbb P^n$.   $R \in \mathbb P^n(\mathbb C) $ does not intersect the set of points  $\{(1,0,\hdots,0),\hdots,(0,\hdots,0,1)\} $. 
  If $V$ is of log-general type, then  there exists a proper Zariski closed subset  $Z$ of  $ X$ 
such that the image of every entire curve
$f:\mathbb C\to V$ is contained in $Z$.
\end{theorem}
We note that the  algebraic degeneracy of entire curves into a complex algebraic variety of positive log Kodaira dimension that admits a  finite morphism onto a semi-abelian has been proved by Noguchi, Winkelman and Yamonoi in \cite{noguchi2007degeneracy}.

\bigskip

\noindent\textbf{Acknowledgements.} 
The authors extend their gratitude to Aaron Levin for providing invaluable insights that greatly improved the reduction method in the proof of Theorem \ref{main_thm}, thereby enabling the theorem's completion for $n\ge 3$. Additionally, the authors would like to express their appreciation to Amos Turchet for engaging discussions and valuable suggestions. Lastly, the first and third authors would like to extend their thanks to the Institute of Mathematics at Academia Sinica for their warm hospitality during their post-doc tenure.

\end{document}